\documentclass[3p]{elsarticle}

\usepackage{amssymb}
\usepackage{amsmath}

\usepackage{xcolor}
\usepackage{graphicx} 
\usepackage[colorlinks=true, linkcolor=blue, citecolor=blue]{hyperref}
\usepackage{algorithm}
\usepackage{algorithmic}
\usepackage{subcaption}
\usepackage{booktabs}

\begin{document}

\begin{frontmatter}



\title{Sequential Bayesian Optimal Experimental Design in Infinite Dimensions via Policy Gradient Reinforcement Learning}


\author[1]{Kaichen Shen} 
\ead{kshen77@gatech.edu}

\author[2]{Peng Chen\corref{cor1}}
\ead{pchen402@gatech.edu}

\cortext[cor1]{Corresponding author.}

\affiliation[1]{
organization={School of Mathematics, Georgia Institute of Technology},
addressline={686 Cherry Street},
city={Atlanta},
postcode={30332}, 
state={Georgia},
country={USA}}

\affiliation[2]{
organization={School of Computational Science and Engineering, Georgia Institute of Technology},
addressline={756 West Peachtree Street Northwest},
city={Atlanta},
postcode={30308}, 
state={Georgia},
country={USA}}

\begin{abstract}
Sequential Bayesian optimal experimental design (SBOED) for PDE-governed inverse problems is computationally challenging, especially for infinite-dimensional random field parameters. High-fidelity approaches require repeated forward and adjoint PDE solves inside nested Bayesian inversion and design loops. We formulate SBOED as a finite-horizon Markov decision process and learn an amortized design policy via policy-gradient reinforcement learning (PGRL), enabling online design selection from the experiment history without repeatedly solving an SBOED optimization problem. 
To make policy training and reward evaluation scalable, we combine dual dimension reduction---active subspace projection for the parameter and principal component analysis for the state---with an adjusted derivative-informed latent attention neural operator (LANO) surrogate that predicts both the parameter-to-solution map and its Jacobian. We use a Laplace-based D-optimality reward while noting that, in general, other expected-information-gain utilities such as KL divergence can also be used within the same framework. We further introduce an eigenvalue-based evaluation strategy that uses prior samples as proxies for maximum a posteriori (MAP) points, avoiding repeated MAP solves while retaining accurate information-gain estimates. Numerical experiments on sequential multi-sensor placement for contaminant source tracking demonstrate approximately $100\times$ speedup over high-fidelity finite element methods, improved performance over random sensor placements, and physically interpretable policies that discover an ``upstream'' tracking strategy.
\end{abstract}


\end{frontmatter}





\section{Introduction} \label{sec:introduction}

Mathematical models governed by partial differential equations (PDEs) are central to predicting complex physical systems, yet their parameters are often uncertain and must be inferred from limited data \cite{stuart2010inverse,ghattas2021learning}. In many applications, experiments are expensive, risky, or limited in time and spatial coverage, so data must be gathered strategically rather than exhaustively. Optimal experimental design (OED) provides a principled way to select experiments that maximize information or reduce uncertainty under such constraints. In the Bayesian setting, Bayesian OED (BOED) commonly targets expected information gain (EIG), formalized as an expected Kullback--Leibler (KL) divergence between posterior and prior \cite{chaloner1995bayesian, rainforth2024modern}. In practice, especially for PDE-constrained inverse problems, EIG is often operationalized through tractable approximations such as covariance-based criteria (e.g., A-/D-optimality) under Laplace and low-rank structure \cite{alexanderian2016on, alexanderian2021optimal, wu2023fast}. In this work, we focus on a Laplace-based D-optimality objective to enable scalable sequential design while still capturing the dominant uncertainty reduction in data-informed subspaces.

Despite decades of progress, BOED for large-scale PDE-constrained inverse problems remains computationally challenging. For each candidate design, evaluating the utility can require repeated forward and adjoint PDE solves, Jacobian computations, and (in Laplace settings) repeated maximum a posteriori (MAP) optimizations across many design and data samples \cite{beck2018fast, wu2023offlineonline, wu2023largescale}. To reduce these costs, practical algorithms exploit Laplace approximations and low-rank structure in Hessians to compute design criteria efficiently \cite{alexanderian2016fast, long2013fast, long2015fast, saibaba2017randomized, attia2018goaloriented, beck2020multilevel, long2015laplace}. Complementary approaches further amortize expensive forward evaluations using polynomial chaos and surrogate modeling \cite{huan2013simulationbased, huan2014gradientbased}. Even with these accelerations, utility evaluation and design optimization can remain costly when designs must be queried repeatedly at scale.

Beyond Laplace-based criteria, a complementary line of work estimates mutual information directly via variational bounds or neural ratio/MI estimators, reducing nested Monte Carlo cost for implicit or nonlinear models \cite{foster2019variational, kleinegesse2020bayesian, go2022robust,orozco2024probabilistic}. While these approaches provide important alternatives for KL-based EIG estimation, our focus here is on a Laplace/D-optimality pipeline that admits fast low-rank reward evaluation and integrates naturally with scalable PDE inverse-problem structure.

Design optimization itself can be combinatorial (e.g., selecting sensor subsets) or continuous (e.g., optimizing sensor locations/times under constraints), and is often nonconvex when each design query is expensive. Gradient-based relaxations, greedy selection, and swapping strategies offer practical compromises that are scalable and often achieve near-optimal designs \cite{alexanderian2014aoptimal, huan2014gradientbased, jagalurmohan2021batch, helin2022edgepromoting, krause2012nearoptimal, wu2023offlineonline, wu2023fast, wu2023largescale}. These considerations motivate \emph{amortized} decision rules that can propose designs quickly once trained, rather than repeatedly solving a new optimization problem online.

When experiments are performed sequentially, SBOED adds an additional layer of difficulty: after each experiment, the posterior becomes the prior for the next decision, so the design policy must balance feedback with planning under future uncertainty. Sequential strategies range from greedy or myopic designs to multi-step lookahead formulations that can be cast as finite-horizon decision processes \cite{drovandi2014sequential, huan2016sequential}. Reinforcement learning provides a natural computational tool for such problems, and recent work uses policy-gradient methods to learn adaptive design policies from simulated episodes \cite{foster2021deep, ivanova2021implicit, blau2022optimizing, shen2025variational, shen2023bayesian}. However, training such policies typically requires many environment interactions, making repeated high-fidelity PDE solves impractical; this motivates surrogate models that can provide fast state and sensitivity evaluations during policy optimization.

Given these cumulative demands, we leverage surrogate modeling based on neural operators. Neural operators provide efficient approximations of PDE solution operators in high dimensions, enabling fast evaluations for new parameter fields and discretizations \cite{kovachki2023neural, li2020fourier, lu2019deeponet, li2024scalable, olearyroseberry2022learning}. For time-dependent PDEs, attention-based architectures improve temporal and multiscale modeling and can capture long-range dependencies \cite{vaswani2017attention, li2022transformer, hao2023gnot, ovadia2024vito}. Crucially, BOED utilities under Laplace approximations---and in particular D-optimality---depend on sensitivities of the parameter-to-observable map; accurate reward evaluation therefore requires \emph{derivative-informed} surrogates trained with Jacobian information \cite{olearyroseberry2024derivativeinformed, olearyroseberry2022derivativeinformed, qiu2024derivativeenhanced, luo2025efficient}. LANO-based surrogates combine dimension reduction, attention, and derivative-informed training to support fast and accurate reward evaluation in SBOED \cite{go2025accurate, go2025sequential}. These developments build on the observation that Jacobians and Hessians in PDE-constrained inverse problems often exhibit low-rank structure concentrated in data-informed subspaces \cite{buithanh2012extremescale, chen2017hessianbased, chen2021bayesian}.

\textbf{Contributions.} We develop a fast and scalable amortized framework for infinite-dimensional SBOED that combines policy-gradient reinforcement learning (PGRL) with derivative-informed surrogate modeling. Our main contributions are:
\begin{itemize}
    \item We cast PDE-constrained SBOED with infinite-dimensional random field parameters as a finite-horizon Markov decision process and train a PGRL policy offline; once trained, it proposes new experiments online without repeatedly solving an SBOED optimization problem.
    \item We adopt D-optimality as the design objective and develop a scalable reward evaluation strategy that exploits low-rank Hessian structure. To further accelerate training, we introduce an approximation that evaluates the spectral properties of the posterior covariance at prior samples---acting as efficient proxies for the MAP point---thereby avoiding expensive iterative MAP optimization inside the RL environment while retaining accurate information-gain estimates
    \item We develop an adjusted derivative-informed latent attention neural operator (LANO) surrogate that predicts both the solution and its Jacobian, coupled with dual dimension reduction (active subspace for parameters and principal component analysis (PCA) for state variables) and additional linear transformations and normalization tailored to the random initial field setting.
    \item We demonstrate the proposed framework on a contaminant source tracking problem with sequential optimal sensor placement (multiple sensors and observations per stage), achieving substantial speedups over high-fidelity PDE solves and learning physically interpretable placement policies.
\end{itemize}

Compared with prior work, our framework is closely related to \cite{go2025sequential, shen2023bayesian} but differs in several key ways:
\begin{itemize}
    \item We target PDE-constrained Bayesian inverse problems with infinite-dimensional random field parameters and richer observation/design settings (e.g., multiple sensors/design variables per stage), while \cite{shen2023bayesian} focuses on a low-dimensional parameter setting with a single observation and design variable per stage and \cite{go2025sequential} focuses on designing optimal observation times.
    \item We learn a PGRL policy offline and propose new experiments online without repeatedly solving an SBOED optimization problem, as in \cite{shen2023bayesian}, whereas \cite{go2025sequential} designs experiments online by solving optimization problems accelerated by a surrogate.
    \item While KL-based EIG is a common utility as in \cite{go2025sequential} and \cite{shen2023bayesian}, we adopt D-optimality and accelerate reward evaluation via an eigenvalue-based strategy that uses prior samples as proxies for MAP points, thus avoiding repeated MAP solves in the design loop that are needed in \cite{go2025sequential}.
    \item Compared to \cite{shen2023bayesian}'s parameter-to-state neural network surrogate for low-dimensional problems, we employ LANO, a derivative-informed neural operator capable of handling infinite-dimensional fields and outputting both states and Jacobians, which is critical for scalable D-optimality computation.
\end{itemize}

The remainder of the paper is organized as follows. Section~\ref{section: Problem formulation} introduces the SBOED setting and MDP formulation, Section~\ref{section: Policy gradient reinforcement learning} presents the PGRL method, Section~\ref{section: Scalable approximations for SBOED} develops the Laplace and low-rank approximations for the reward, and Section~\ref{section: Surrogate modeling} describes the LANO surrogate and training. We then present the contaminant source tracking setup in Section~\ref{section: Case study: contaminant source tracking}, report numerical experiments in section~\ref{section: Numerical experiments}, and conclude in Section~\ref{section: conclusion}.


\section{Problem formulation}
\label{section: Problem formulation}

We consider a time-dependent PDE with an uncertain random field parameter $m$ and seek to infer $m$ from noisy observations of the state $u(t,\cdot)$. At each stage $n=1,\dots,N$, we choose a design vector $\mathbf{d}_n$ (e.g., sensor placements) that induces an observation operator, collect the resulting data $\mathbf{y}_n$, and update the posterior distribution of $m$. The SBOED objective is to choose $\{\mathbf{d}_n\}_{n=1}^N$ (or equivalently a policy) to maximize cumulative information gain, which we quantify using D-optimality under a Laplace approximation.

\subsection{Infinite-dimensional Bayesian inverse problem}
\label{subsection: Infinite-dimensional Bayesian inverse problem}

We consider a system governed by a time-dependent PDE with an infinite-dimensional uncertain parameter:
\begin{align}
    \label{eq: pde}
    \partial_{t} u(t, x) + \mathcal{R}(u(t, x), m(x)) = 0, \quad (t, x) \in (0, T] \times \Omega,
\end{align}
where $T > 0$ denotes the terminal time and $\Omega \subset \mathbb{R}^{d_x}$ is an open, bounded physical domain. The state variable $u(t, \cdot)$ lies in a Hilbert space $\mathcal{V}$ defined on $\Omega$, satisfying appropriate boundary conditions for all $t \in (0, T)$ and the initial condition $u(0, \cdot) = u_{0}$, which may depend on the parameter $m$ as in our case study. The parameter $m \in \mathcal{M}$ is a random field in the Hilbert space $\mathcal{M}$ that may represent distributed sources, initial conditions, boundary conditions, coefficients, or geometry. The operator $\mathcal{R}: \mathcal{V} \times \mathcal{M} \rightarrow \mathcal{V}'$ denotes a differential operator, with $\mathcal{V}'$ representing the dual space of $\mathcal{V}$. 

Let $0 < t_{1} < \cdots < t_{K} \leq T$ denote the $K$ observation times at which data are collected; these times need not include $t=0$. For notational convenience, we write $u_{k} = u(t_{k}, \cdot)$ for $k = 1, \cdots, K$, and let $u_{k}(m)$ denote the solution to \eqref{eq: pde} at time $t_{k}$ given the parameter $m$.
We define the parameter-to-observable (PtO) map at time $t_{k}$ as $\mathcal{F}_{k}: \mathcal{M} \rightarrow \mathbb{R}^{d_y}$ for $k = 1, \cdots, K$. This is typically expressed as $\mathcal{F}_{k}(m) = \mathcal{B}_{k}(u_{k}(m))$, where $\mathcal{B}_{k}: \mathcal{V} \rightarrow \mathbb{R}^{d_y}$ is the observation operator at time $t_{k}$ (depending on design variables). 
In the sequential design setting, the stage-$n$ design vector $\mathbf{d}_n$ parameterizes the observation operator(s) used at that stage (e.g., sensor locations), which we denote by $\mathcal{B}_n$ (and apply at the corresponding observation indices $k \in \mathcal{K}_n$).
The observation vector $\mathbf{y}^{(k)}$ at time $t_{k}$ is modeled with additive noise:
\begin{align}
\label{eq: observation}
\mathbf{y}^{(k)} = \mathcal{F}_{k}(m) + \boldsymbol{\epsilon}^{(k)},
\end{align}
where each noise term $\boldsymbol{\epsilon}^{(k)}$ follows a centered Gaussian distribution $\mathcal{N}(\mathbf{0}, \boldsymbol{\Gamma}_{\text{noise}})$ with covariance matrix $\boldsymbol{\Gamma}_{\text{noise}} \in \mathbb{R}^{d_y \times d_y}$. Assuming statistical independence across time steps $t_{k}$, the likelihood of the aggregated observation data $\mathbf{y} = (\mathbf{y}^{(1)}, \cdots, \mathbf{y}^{(K)})$ is given by:
\begin{align}
\label{eq: likelihood}
p_{\text{like}}(\mathbf{y} | m) \propto \exp \left( - \frac{1}{2} \sum_{k=1}^{K} \| \mathbf{y}^{(k)} - \mathcal{F}_{k}(m) \|^2_{\boldsymbol{\Gamma}_{\text{noise}}^{-1}} \right).
\end{align}

The uncertain parameter $m$ is endowed with a Gaussian prior measure:
\begin{align}
\label{eq: prior}
\mu_{\text{prior}} = \mathcal{N}(m_{\text{prior}}, \mathcal{C}_{\text{prior}}),
\end{align}
where $m_{\text{prior}}$ is the mean field defined on $\Omega$. The covariance operator is of the Matérn class, defined as $\mathcal{C}_{\text{prior}} = \mathcal{A}^{-\alpha}$, where $\mathcal{A} = -\gamma \Delta + \delta I$ is an elliptic operator with appropriate boundary conditions \cite{lindgren2011explicit}. The parameter $\alpha > d_x / 2$ ensures that $\mathcal{C}_{\text{prior}}$ is of trace class. The hyperparameters $\alpha, \gamma, \delta > 0$ collectively determine the smoothness, variance, and correlation length of the random field.

Given the observation data $\mathbf{y}$, the posterior distribution $\mu_{\text{post}}$ of the random field $m$ is characterized by Bayes' rule through the Radon--Nikodym derivative:
\begin{align}
\label{eq: Bayes' rule}
\frac{d\mu_{\text{post}}}{d\mu_{\text{prior}}}(m) = \frac{1}{p(\mathbf{y})} p_{\text{like}}(\mathbf{y} | m),
\end{align}
where $p(\mathbf{y})$ is the normalization constant (or model evidence):
\begin{align}
\label{eq: evidence}
p(\mathbf{y}) = \int_{\mathcal{M}} p_{\text{like}}(\mathbf{y} | m) d\mu_{\text{prior}}(m),
\end{align}
which is typically computationally intractable to evaluate. The primary objective of the Bayesian inverse problem is to characterize and sample from the posterior $\mu_{\text{post}}$ conditioned on the observations $\mathbf{y}$.

\subsection{Sequential Bayesian optimal experimental design} \label{subsection: Sequential Bayesian optimal experimental design}

We now formalize the SBOED problem within the Bayesian inverse framework of Section~\ref{subsection: Infinite-dimensional Bayesian inverse problem}. Consider a sequential experimental process with $N$ stages, indexed by $n = 1, \dots, N$. At stage $n$, a design vector $\mathbf{d}_n \in \mathcal{D}_n \subseteq \mathbb{R}^{d_d}$ is selected. Each experiment may include multiple observation times; let $\mathcal{K}_n$ denote the set of (consecutive) indices of those times with $|\mathcal{K}_n| \geq 1$ (e.g., $\mathcal{K}_1 = \{1,2,3\}$ and $\mathcal{K}_2 = \{4,5\}$). We write $\mathbf{y}_n = \{\mathbf{y}^{(k)}\}_{k \in \mathcal{K}_n}$ for the corresponding data and $K = \sum_{n=1}^N |\mathcal{K}_n|$ for the total number of observation times. 

Let $\mathbf{d}_{1:n} = (\mathbf{d}_1, \dots, \mathbf{d}_n)$ and $\mathbf{y}_{1:n} = (\mathbf{y}_1, \dots, \mathbf{y}_n)$ denote the accumulated histories of design and observation vectors, respectively, up to the $n$-th stage. The posterior distribution of the random field parameter $m$, conditioned on these histories, is denoted by $\mu(m | \mathbf{d}_{1:n}, \mathbf{y}_{1:n})$ and is obtained via Bayes' rule \eqref{eq: Bayes' rule}.

In general, the information gain in Bayesian experimental design is often quantified by the expected KL divergence between posterior and prior (expected information gain). In this work, we adopt the D-optimality criterion \cite{alexanderian2016on, wu2023fast}, which measures the contraction of the posterior covariance. The incremental information gain from the $n$-th experiment is quantified by applying D-optimality to successive posteriors:
\begin{align}
\label{eq: D-optimality}
D_{\text{opt}}(\mu(m | \mathbf{d}_{1:n}, \mathbf{y}_{1:n}) \parallel \mu(m | \mathbf{d}_{1:n-1}, \mathbf{y}_{1:n-1})).
\end{align}
Here $\mu(m | \mathbf{d}_{1:0}, \mathbf{y}_{1:0})$ represents the initial prior $\mu_{\text{prior}}$ defined in \eqref{eq: prior}. A detailed derivation of the D-optimality expression under Laplace approximation of the posterior is provided in Section~\ref{subsection: Computation of the reward function}.

The objective of the SBOED problem is to determine the optimal sequence of design vectors $\mathbf{d}^*_{1:N}$ that maximizes the expected cumulative information gain over all $N$ experiments. This optimization problem is formally stated as:
\begin{align}
\label{eq: SBOED goal}
\max_{\mathbf{d}_{1:N}} \mathbb{E}_{\mathbf{y}_{1:N} | \mathbf{d}_{1:N}} \left[ \sum_{n=1}^{N} D_{\text{opt}} \left( \mu(m | \mathbf{d}_{1:n}, \mathbf{y}_{1:n}) \parallel \mu(m | \mathbf{d}_{1:n-1}, \mathbf{y}_{1:n-1}) \right) \right].
\end{align}
The stochasticity of the observation sequence $\mathbf{y}_{1:N}$ in \eqref{eq: SBOED goal} is induced by the uncertainty in the random field $m$ and the observation noise $\boldsymbol{\epsilon}_{1:N}$. 
To approximate the expectation in \eqref{eq: SBOED goal}, we generate synthetic observations following \cite{shen2023bayesian}: first sample $m \sim \mu_{\text{prior}}$ in \eqref{eq: prior}, then solve \eqref{eq: pde} to obtain the state at the required time points $t_k$, and finally apply the observation operators induced by the designs $\mathbf{d}_{1:N}$ and add noise as in \eqref{eq: observation}.

\subsection{Markov decision process formulation of SBOED} \label{subsection: Markov decision process formulation of SBOED}

The objective in \eqref{eq: SBOED goal} provides a general formulation wherein design vectors $\mathbf{d}_{1:N}$ are selected from an arbitrary space. To make the sequential structure explicit, we cast the SBOED framework as a Markov decision process (MDP) \cite{shen2023bayesian}. Formally, this MDP is defined by the tuple $(\mathcal{X}, \mathcal{D}, \mathcal{T}, g)$, representing the state space, design space, transition dynamics, and reward function, respectively. A policy is a mapping from the state space to the design space. We detail these components below.

\textbf{State.} The state vector $\mathbf{x}_{n} = [\mathbf{x}_{n,b}, \mathbf{x}_{n, p}] \in \mathcal{X}_{n}$ characterizes the system status prior to the $n$-th experiment. This vector consists of the belief state $\mathbf{x}_{n, b}$, which encapsulates the uncertainty associated with the random field parameter $m$, and the physical state $\mathbf{x}_{n, p}$, which contains all remaining deterministic information. Consistent with \cite{shen2023bayesian}, we define the state vector $\mathbf{x}_{n}$ as the concatenation of the history of design and observation vectors preceding the $n$-th experiment:
\begin{align}
\label{eq: state variable}
& \mathbf{x}_{n} = [\mathbf{d}_{1}, \mathbf{y}_{1}, \cdots, \mathbf{d}_{n-1}, \mathbf{y}_{n-1}] \in \mathcal{X}_{n},
\end{align}
where $\mathbf{x}_{1}$ is initialized as an empty vector. This concatenation is formally referred to as the information vector and constitutes a sufficient statistic for the posterior distribution $\mu(m | \mathbf{d}_{1:n-1}, \mathbf{y}_{1:n-1})$.

\textbf{Design and policy.} The policy $\pi_{n}$ governing the $n$-th experiment is defined as a deterministic mapping from the state space to the design space:
\begin{align}
\label{eq: policy nth experiment}
\pi_{n}: \mathcal{X}_{n} \rightarrow \mathcal{D}_{n}: \mathbf{x}_{n} \mapsto \mathbf{d}_{n}.
\end{align}
We denote the complete sequence of policies as $\pi = \{ \pi_{1}, \cdots, \pi_{N} \}$.

\textbf{Transition dynamics.} The transition function $\mathcal{T}$ governs the evolution from state $\mathbf{x}_{n}$ to $\mathbf{x}_{n+1}$ following the execution of design $\mathbf{d}_{n}$ and the acquisition of observation $\mathbf{y}_{n}$. In accordance with \eqref{eq: state variable}, the transition dynamics are modeled by appending the current design-observation pair to the existing state vector:
\begin{align}
\label{eq: transition dynamics}
& \mathbf{x}_{n + 1} = \mathcal{T}(\mathbf{x}_{n}, \mathbf{d}_{n}, \mathbf{y}_{n}) := [\mathbf{x}_{n}, \mathbf{d}_{n}, \mathbf{y}_{n}] \in \mathcal{X}_{n + 1}.
\end{align}

\textbf{Reward function.} 
We define the immediate reward $g_{n}(\mathbf{x}_{n}, \mathbf{d}_{n}, \mathbf{y}_{n}) \in \mathbb{R}$ as the utility derived from conducting the $n$-th experiment with design $\mathbf{d}_{n}$ and observing $\mathbf{y}_{n}$ from state $\mathbf{x}_{n}$. Additionally, let $g_{N + 1}(\mathbf{x}_{N+1}) \in \mathbb{R}$ represent the terminal reward evaluated upon completion of the experimental sequence, where $\mathbf{x}_{N+1}$ aggregates the full history of designs and observations:
\begin{align}
    \label{eq: terminal reward}
    \mathbf{x}_{N+1} = [\mathbf{d}_{1}, \mathbf{y}_{1}, \cdots, \mathbf{d}_{N}, \mathbf{y}_{N}].
\end{align}
In the context of the SBOED problem for Bayesian inverse problems defined in \eqref{eq: SBOED goal}, the reward function is quantified by D-optimality. We specifically define the \textit{incremental formulation} of the reward as:
\begin{align}
    \label{eq: incremental formulation immediate}
   g_{n}(\mathbf{x}_{n}, \mathbf{d}_{n}, \mathbf{y}_{n}) &:= D_{\text{opt}}(\mu(m | \mathbf{x}_{n+1}) \parallel \mu(m | \mathbf{x}_{n})), \\
    \label{eq: incremental formulation terminal}
    g_{N+1}(\mathbf{x}_{N+1}) &:= 0,
\end{align}
which performs Bayesian inference after each experiment. Alternatively, we consider the \textit{terminal formulation} of the reward, which assigns utility solely at the end of the process:
\begin{align}
    \label{eq: terminal formulation immediate}
    g_{n}(\mathbf{x}_{n}, \mathbf{d}_{n}, \mathbf{y}_{n}) &:=0, \\
    \label{eq: terminal formulation terminal}
    g_{N+1}(\mathbf{x}_{N+1}) &:= D_{\text{opt}}(\mu(m | \mathbf{x}_{N+1}) \parallel \mu(m | \mathbf{x}_{1})).
\end{align}
The terminal formulation offers significant computational advantages by requiring only a single Bayesian inference operation upon completion of the sequence.

\textbf{MDP formulation of SBOED.} 
Consequently, the objective of the SBOED problem is to identify the optimal policy sequence $\pi$ that maximizes the expected cumulative reward. This is formally stated as the following optimization problem:
\begin{align}
    \label{eq: SBOED MDP formulation}
    & \max_{\pi = \{ \pi_{1}, \cdots, \pi_{N} \}} U(\pi) := \mathbb{E}_{\mathbf{y}_{1:N} | \pi} \left[ \sum_{n=1}^{N} g_{n}(\mathbf{x}_{n}, \mathbf{d}_{n}, \mathbf{y}_{n}) + g_{N+1}(\mathbf{x}_{N+1}) \right], \nonumber \\
    & \text{s.t. } \mathbf{x}_{1}  = \emptyset, \; \mathbf{d}_{n} = \pi_{n}(\mathbf{x}_{n}), \;  \mathbf{x}_{n+1} = \mathcal{T}(\mathbf{x}_{n}, \mathbf{d}_{n}, \; \mathbf{y}_{n}), n = 1, \cdots, N.
\end{align}
This MDP formulation explicitly realizes the general objective in \eqref{eq: SBOED goal} by constraining the design vectors $\mathbf{d}_{1:N}$ to be outputs of the policies $\pi$. This structure facilitates the application of policy gradient reinforcement learning methods (Section~\ref{section: Policy gradient reinforcement learning}). We demonstrate in Section~\ref{subsection: Computation of the reward function} that the incremental and terminal reward formulations are equivalent for D-optimality.


\section{Policy gradient reinforcement learning} 
\label{section: Policy gradient reinforcement learning}

This section applies policy gradient reinforcement learning (PGRL) to the SBOED problem cast as the Markov decision process in \eqref{eq: SBOED MDP formulation}. By parameterizing the policies, we obtain a finite-dimensional optimization problem that can be solved using gradient ascent. To evaluate the policy gradient, we introduce a Q-function (critic) that estimates expected future rewards.

\textbf{Episode view.} One training episode proceeds as follows: sample a parameter realization (and observation noise), roll out $N$ sequential designs using the current policy to generate observations, compute the terminal D-optimality reward from the resulting history, and then update the policy (actor) and Q-function (critic) using stochastic gradients. In our implementation, we use the terminal reward formulation in \eqref{eq: terminal formulation immediate}--\eqref{eq: terminal formulation terminal} to avoid repeated Bayesian inversions during rollout.

\subsection{Derivation of the policy gradient} 
\label{subsection: Derivation of the policy gradient}

We parameterize the stage-$n$ policy $\pi_n$ by $\mathbf{w}_n$ and denote it by $\pi_{\mathbf{w}_n}$. Collecting all parameters yields $\mathbf{w} = [\mathbf{w}_1, \dots, \mathbf{w}_N]^\top$ and the corresponding policy sequence $\pi_{\mathbf{w}}$. The MDP formulation \eqref{eq: SBOED MDP formulation} is thus restated as a maximization problem over $\mathbf{w}$:
\begin{align}
    \label{eq: parameterized SBOED MDP formulation}
    & \max_{\mathbf{w} = [\mathbf{w}_1, \dots, \mathbf{w}_N]} U(\mathbf{w}) := \mathbb{E}_{\mathbf{y}_{1:N} | \pi_{\mathbf{w}}} \left[ \sum_{n=1}^{N} g_{n}(\mathbf{x}_{n}, \mathbf{d}_{n}, \mathbf{y}_{n}) + g_{N+1}(\mathbf{x}_{N+1}) \right], \nonumber \\
    & \text{s.t.} \; \mathbf{x}_{1} = \emptyset, \; \mathbf{d}_{n} = \pi_{\mathbf{w}_n}(\mathbf{x}_{n}), \; \mathbf{x}_{n+1} = \mathcal{T}(\mathbf{x}_{n}, \mathbf{d}_{n}, \mathbf{y}_{n}), \; n = 1, \dots, N.
\end{align}
We maximize $U(\mathbf{w})$ via gradient ascent. As derived in \cite{shen2023bayesian}, the gradient of the objective with respect to the parameters $\mathbf{w}$ is given by:
\begin{align}
    \label{eq: policy gradient}
    \nabla_{\mathbf{w}} U(\mathbf{w}) = \mathbb{E}_{\mathbf{y}_{1:N} | \pi_{\mathbf{w}}} \left[ \sum_{n=1}^{N} \nabla_{\mathbf{w}} \pi_{\mathbf{w}_n}(\mathbf{x}_{n}) \nabla_{\mathbf{d}_{n}} Q_{n}^{\pi_{\mathbf{w}}}(\mathbf{x}_{n}, \mathbf{d}_{n}) \right],
\end{align}
where $\nabla_{\mathbf{w}} \pi_{\mathbf{w}_n}(\mathbf{x}_{n})$ is the Jacobian of the policy and $\nabla_{\mathbf{d}_{n}} Q_{n}^{\pi_{\mathbf{w}}}(\mathbf{x}_{n}, \mathbf{d}_{n})$ is the gradient of the Q-function with respect to the design vector. The Q-function $Q_{n}^{\pi_{\mathbf{w}}}$ is defined as the expected remaining cumulative reward when taking design $\mathbf{d}_n$ in state $\mathbf{x}_n$ and following $\pi_{\mathbf{w}}$ thereafter:
\begin{align}
    \label{eq: Q-function}
    Q_{n}^{\pi_{\mathbf{w}}}(\mathbf{x}_{n}, \mathbf{d}_{n}) := \mathbb{E}_{\mathbf{y}_{n:N} | \pi_{\mathbf{w}}, \mathbf{x}_{n}, \mathbf{d}_{n}} \left[ \sum_{\ell=n}^{N} g_{\ell}(\mathbf{x}_{\ell}, \mathbf{d}_{\ell}, \mathbf{y}_{\ell}) + g_{N+1}(\mathbf{x}_{N+1}) \right].
\end{align}
The Q-function satisfies the Bellman recursion:
\begin{align}
    \label{eq: Q-function recursion}
    Q_{n}^{\pi_{\mathbf{w}}}(\mathbf{x}_{n}, \mathbf{d}_{n}) = \mathbb{E}_{\mathbf{y}_{n} | \pi_{\mathbf{w}}, \mathbf{x}_{n}, \mathbf{d}_{n}} \left[ g_{n}(\mathbf{x}_{n}, \mathbf{d}_{n}, \mathbf{y}_{n}) + Q_{n+1}^{\pi_{\mathbf{w}}}(\mathbf{x}_{n+1}, \mathbf{d}_{n+1}) \right],
\end{align}
for $n = 1, \dots, N-1$.

For numerical implementation, the expectation in \eqref{eq: policy gradient} is approximated via Monte Carlo (MC) estimation. Under the current policy $\pi_{\mathbf{w}}$, we generate $M$ independent trajectories $\{ \mathbf{d}_{1:N}^{(i)}, \mathbf{y}_{1:N}^{(i)} \}_{i=1}^M$ by simulating the forward model and measurement process described in Section~\ref{section: Problem formulation}. The gradient is then estimated as:
\begin{align}
    \label{eq: MC policy gradient}
    \nabla_{\mathbf{w}} U(\mathbf{w}) \approx \frac{1}{M} \sum_{i=1}^{M} \sum_{n=1}^{N} \nabla_{\mathbf{w}} \pi_{\mathbf{w}_n}(\mathbf{x}_{n}^{(i)}) \nabla_{\mathbf{d}_{n}} Q_{n}^{\pi_{\mathbf{w}}}(\mathbf{x}_{n}^{(i)}, \mathbf{d}_{n}^{(i)}).
\end{align}

\subsection{Actor-critic method} 
\label{subsection: actor-critic}

Evaluating \eqref{eq: MC policy gradient} requires $\nabla_{\mathbf{w}} \pi_{\mathbf{w}_n}$ and $\nabla_{\mathbf{d}_n} Q_{n}^{\pi_{\mathbf{w}}}$. The policy Jacobian is obtained via automatic differentiation of the parameterized policy network. To compute $\nabla_{\mathbf{d}_n} Q$, we adopt an actor-critic architecture and approximate the Q-function with a differentiable critic network, enabling backpropagation through $Q$. We detail the actor (policy) and critic (Q-function) parameterizations below.

\textbf{Policy network.} 
We approximate the sequence of policies $\pi_{\mathbf{w}} = \{ \pi_{\mathbf{w}_1}, \dots, \pi_{\mathbf{w}_N} \}$ using a unified neural network architecture, denoted as $\mathcal{P}_{\mathbf{w}}(n, \mathbf{x}_n)$. We utilize a single network for all stages, treating the stage index $n$ as an explicit input. The network inputs comprise the one-hot encoded stage index $n$ and the current state vector $\mathbf{x}_n$, defined in \eqref{eq: state variable} as the history of designs and observations.

Since the dimension of $\mathbf{x}_n$ grows linearly with $n$, reaching a maximum of $(N-1)(d_d + d_y)$ at $n=N$, we employ zero-padding to standardize the input size. The state vector is padded to the maximum dimensionality to construct a fixed-size input vector:
\begin{align}
    \label{eq: state variable neural network input}
    [\mathbf{d}_1, \dots, \mathbf{d}_{n-1}, \mathbf{0}, \dots, \mathbf{0}, \mathbf{y}_1, \dots, \mathbf{y}_{n-1}, \mathbf{0}, \dots, \mathbf{0}] \in \mathbb{R}^{(N-1)(d_d + d_y)}.
\end{align}
This padded information vector is a sufficient statistic for the posterior under the assumed model \cite{shen2023bayesian}, and thus yields a Markov state representation for the sequential design process.
The policy parameters $\mathbf{w}$ are updated via the gradient ascent direction in \eqref{eq: MC policy gradient}.

\textbf{Q-network.} 
The Q-function is similarly parameterized by a neural network (critic), denoted by $\mathcal{Q}_{\boldsymbol{\eta}}(n, \mathbf{x}_n, \mathbf{d}_n)$, where $\boldsymbol{\eta}$ represents the trainable parameters. In contrast to the policy network, the critic takes the design vector $\mathbf{d}_n$ as an additional input appended to the padded state vector in \eqref{eq: state variable neural network input}. 

The Q-network is trained via supervised learning to satisfy the Bellman recursion. We define the loss function $\mathcal{L}(\boldsymbol{\eta})$ based on the mean squared error (MSE) between the predicted Q-values and the temporal difference (TD) targets derived from \eqref{eq: Q-function recursion}:
\begin{align}
    \label{eq: Q-network training loss function}
    \mathcal{L}(\boldsymbol{\eta}) := \frac{1}{M} \sum_{i=1}^{M} \sum_{n=1}^{N-1} \left( \mathcal{Q}_{\boldsymbol{\eta}}^{(n,i)} - \left(g_n^{(i)} + \mathcal{Q}_{\boldsymbol{\eta}}^{(n+1, i)} \right) \right)^2 + \frac{1}{M} \sum_{i=1}^{M} \left( \mathcal{Q}_{\boldsymbol{\eta}}^{(N,i)} - \left( g_N^{(i)} + g_{N+1}^{(i)} \right) \right)^2,
\end{align}
where $\mathcal{Q}_{\boldsymbol{\eta}}^{(n,i)} = \mathcal{Q}_{\boldsymbol{\eta}}(n, \mathbf{x}_n^{(i)}, \mathbf{d}_n^{(i)})$ and $g_n^{(i)} = g_n(\mathbf{x}_n^{(i)}, \mathbf{d}_n^{(i)}, \mathbf{y}_n^{(i)})$, with $i = 1, \dots, M$ and $n = 1, \dots, N-1$.
The samples $\{ \mathbf{d}_{1:N}^{(i)}, \mathbf{y}_{1:N}^{(i)} \}_{i=1}^M$ used in this loss function are the same trajectories generated for the policy gradient estimation, ensuring data efficiency.

\section{Scalable approximations for SBOED}
\label{section: Scalable approximations for SBOED}

This section presents scalable approximations for solving the Bayesian inverse problem \eqref{eq: Bayes' rule} and evaluating the D-optimality reward in \eqref{eq: incremental formulation immediate} and \eqref{eq: terminal formulation terminal}. These ingredients enable us to extend the framework of \cite{shen2023bayesian}, which is restricted to low-dimensional parameter spaces, to settings with infinite-dimensional random field parameters. We approximate three core components required for scalable reward evaluation: (i) the posterior via a Laplace approximation (MAP point and local curvature), (ii) the dominant generalized eigenpairs of the Gauss--Newton Hessian via low-rank structure, and (iii) the D-optimality reward via the resulting eigenvalues.
Unless otherwise noted, calligraphic symbols denote function-space objects/operators (e.g., $\mathcal{M}$, $\mathcal{C}_{\text{prior}}$), while bold symbols denote their finite-dimensional discretizations (e.g., $\mathbf{m}$, $\mathbf{C}_{\text{prior}}$).

\subsection{High-fidelity discretization} \label{subsection: High-fidelity discretization}

We discretize the infinite-dimensional random field $m \in \mathcal{M}$ using the finite element method (FEM) in a finite-dimensional subspace $\mathcal{M}_{d_m} \subset \mathcal{M}$ of dimension $d_m$. Let $\{\phi_j\}_{j=1}^{d_m}$ denote the piecewise continuous Lagrange basis functions spanning $\mathcal{M}_{d_m}$. For a mesh on $\Omega$ with vertices $\{v_j\}_{j=1}^{d_m}$, these basis functions satisfy $\phi_j(v_{j'}) = \delta_{jj'}$. The approximation $\hat{m}$ is expressed as:
\begin{align}
    \label{eq: parameter discretization}
    \hat{m}(x) = \sum_{j=1}^{d_m} m_j \phi_j(x).
\end{align}
We denote the coefficient vector as $\mathbf{m} = (m_1, \dots, m_{d_m})^\top \in \mathbb{R}^{d_m}$. To ensure a high-fidelity representation of the random field, the dimension $d_m$ is typically large.

The Gaussian prior $\mu_{\text{prior}}$ in \eqref{eq: prior} is discretized into a multivariate normal distribution $\mathcal{N}(\mathbf{m}_{\text{prior}}, \mathbf{C}_{\text{prior}})$, where $\mathbf{m}_{\text{prior}}$ contains the nodal values of the mean field. For the case $\alpha=2$, the Matérn covariance discretization takes the form $\mathbf{C}_{\text{prior}} = \mathbf{A}^{-1} \mathbf{M} \mathbf{A}^{-1}$, where $\mathbf{M} \in \mathbb{R}^{d_m \times d_m}$ is the mass matrix and $\mathbf{A} \in \mathbb{R}^{d_m \times d_m}$ is the stiffness-mass discretization of the associated elliptic operator.

We similarly discretize the state variables (potentially in a different finite element space). Let $\mathbf{u}_k(\mathbf{m}) \in \mathbb{R}^{d_u}$ denote the discrete state vector at time $t_k$, and let $\mathbf{B}_k: \mathbb{R}^{d_u} \to \mathbb{R}^{d_y}$ and $\mathbf{F}_k: \mathbb{R}^{d_m} \to \mathbb{R}^{d_y}$ denote the discretized observation operator and PtO map, respectively.

\subsection{Laplace approximation of the posterior distribution} \label{subsection: Laplace approximation of the posterior distribution}

We approximate the posterior distribution of the discretized parameter $\mathbf{m}$, conditioned on the history $\mathbf{x}_{n+1}$ (i.e., designs $\mathbf{d}_{1:n}$ and observations $\mathbf{y}_{1:n}$), using the Laplace method, yielding a Gaussian approximation:
\begin{align}
    \label{eq: Laplace approximation}
    \mu(\mathbf{m} | \mathbf{x}_{n+1}) \approx \mathcal{N}(\mathbf{m}_{\text{MAP}}^{\mathbf{x}_{n+1}}, \mathbf{C}_{\text{post}}^{\mathbf{x}_{n+1}}), \quad n = 1, \dots, N.
\end{align}
The maximum a posteriori (MAP) estimate, $\mathbf{m}_{\text{MAP}}^{\mathbf{x}_{n+1}}$, is obtained by solving the deterministic optimization problem: 
\begin{align}
    \label{eq: MAP optimization}
    \min_{\mathbf{m} \in \mathbb{R}^{d_m}} \mathcal{J}(\mathbf{m}) := \frac{1}{2} \sum_{\ell=1}^{n} \sum_{k \in \mathcal{K}_\ell} \| \mathbf{y}^{(k)} - \mathbf{F}_k(\mathbf{m}) \|^2_{\boldsymbol{\Gamma}_{\text{noise}}^{-1}} + \frac{1}{2} \| \mathbf{m} - \mathbf{m}_{\text{prior}} \|^2_{\mathbf{C}_{\text{prior}}^{-1}},
\end{align}
where $\mathcal{K}_\ell$ denotes the set of time indices associated with the observations in $\mathbf{y}_\ell$. The objective comprises a data-misfit term weighted by the noise covariance and a quadratic regularization term induced by the prior. We solve \eqref{eq: MAP optimization} efficiently using an inexact Newton-CG algorithm \cite{villa2021hippylib}.

The posterior covariance matrix $\mathbf{C}_{\text{post}}^{\mathbf{x}_{n+1}}$ is given by the inverse of the Hessian of the negative log-posterior evaluated at the MAP point:
\begin{align}
    \label{eq: posterior covariance exact}
    \mathbf{C}_{\text{post}}^{\mathbf{x}_{n+1}} = \left( \mathbf{H}_{\text{misfit}}^{\mathbf{x}_{n+1}}(\mathbf{m}_{\text{MAP}}^{\mathbf{x}_{n+1}}) + \mathbf{C}_{\text{prior}}^{-1} \right)^{-1}.
\end{align}
To accelerate computation, we approximate the Hessian of the misfit term $\mathbf{H}_{\text{misfit}}^{\mathbf{x}_{n+1}} (\mathbf{m}_{\text{MAP}}^{\mathbf{x}_{n+1}})$ at $\mathbf{m}_{\text{MAP}}^{\mathbf{x}_{n+1}}$ using the Gauss-Newton (GN) approximation \cite{wu2023fast}, which discards second-order derivative terms of the PtO map:
\begin{align}
    \label{eq: GN approximation}
    \mathbf{H}_{\text{misfit}}^{\mathbf{x}_{n+1}} (\mathbf{m}_{\text{MAP}}^{\mathbf{x}_{n+1}}) \approx \mathbf{H}_{\text{GN}}^{\mathbf{x}_{n+1}} (\mathbf{m}_{\text{MAP}}^{\mathbf{x}_{n+1}}) = \sum_{\ell=1}^{n} \sum_{k \in \mathcal{K}_\ell} \left( \nabla_{\mathbf{m}} \mathbf{F}_k(\mathbf{m}_{\text{MAP}}^{\mathbf{x}_{n+1}}) \right)^\top \boldsymbol{\Gamma}_{\text{noise}}^{-1} \left( \nabla_{\mathbf{m}} \mathbf{F}_k(\mathbf{m}_{\text{MAP}}^{\mathbf{x}_{n+1}}) \right).
\end{align}

Given the high dimensionality $d_m$, explicitly forming or inverting the Hessian is computationally prohibitive. We therefore employ a low-rank approximation by solving the generalized eigenvalue problem:
\begin{align}
    \label{eq: generalized eigenvalue}
    \mathbf{H}_{\text{GN}}^{\mathbf{x}_{n+1}}(\mathbf{m}_{\text{MAP}}^{\mathbf{x}_{n+1}}) \mathbf{w}_j = \lambda_j \mathbf{C}_{\text{prior}}^{-1} \mathbf{w}_j, \quad j = 1, \dots, r,
\end{align}
using a scalable randomized algorithm \cite{villa2021hippylib}. The eigenvalues are ordered such that $\lambda_1 \ge \dots \ge \lambda_r$, and the spectrum is truncated at rank $r$ where $\lambda_r \ll 1$. The posterior covariance is then approximated as:
\begin{align}
    \label{eq: low rank covariance}
    \mathbf{C}_{\text{post}}^{\mathbf{x}_{n+1}} \approx \mathbf{C}_{\text{prior}} - \mathbf{W}_r \mathbf{D}_r \mathbf{W}_r^\top,
\end{align}
where $\mathbf{W}_r = [\mathbf{w}_1, \dots, \mathbf{w}_r] \in \mathbb{R}^{d_m \times r}$ contains the generalized eigenvectors, and $\mathbf{D}_r = \text{diag}(\frac{\lambda_1}{1+\lambda_1}, \dots, \frac{\lambda_r}{1+\lambda_r})$.

\subsection{Computation of the reward function} \label{subsection: Computation of the reward function}

The D-optimality reward can be approximated by the dominant generalized eigenvalues of the Gauss--Newton Hessian and is therefore well suited to low-rank approximations. In particular, D-optimality (the reduction in the determinant of the covariance, representing the volume contraction of the confidence ellipsoid) is approximated by \cite{alexanderian2016on}:
\begin{align}
    \label{eq: D-optimality computation}
    D_{\text{opt}}(\mu(\mathbf{m} | \mathbf{x}_{n+1}) \parallel \mu(\mathbf{m} | \mathbf{x}_{1})) &= \frac{1}{2} \log \frac{\det(\mathbf{C}_{\text{prior}})}{\det(\mathbf{C}_{\text{post}}^{\mathbf{x}_{n+1}})} \approx \frac{1}{2} \sum_{j=1}^{r} \log(1 + \lambda_j).
\end{align}

As noted in Section~\ref{section: Problem formulation}, the incremental and terminal reward formulations are equivalent for D-optimality. Substituting \eqref{eq: D-optimality computation} into \eqref{eq: SBOED MDP formulation} shows that the incremental logarithmic terms telescope into the terminal log-determinant term. We therefore adopt the terminal formulation to reduce computational cost, since it requires solving the Bayesian inverse problem only once at the end of the experimental sequence.

\textbf{Summary.} Given an experiment history (designs and observations), we compute a Laplace approximation of the posterior by solving for the MAP point and estimating the dominant generalized eigenvalues of the Gauss--Newton Hessian. The D-optimality reward is then evaluated from these eigenvalues using \eqref{eq: D-optimality computation}.

\section{Surrogate modeling}
\label{section: Surrogate modeling}

Although the approximations in Section~\ref{section: Scalable approximations for SBOED} are scalable, they still require many high-fidelity PDE solves, which is prohibitively expensive for training the PGRL policy in Section~\ref{section: Policy gradient reinforcement learning}. To address this cost, we introduce a surrogate model, a latent attention neural operator (LANO), designed to approximate the solution operator of the time-dependent PDE \eqref{eq: pde} and its Jacobian with respect to the parameter. We describe its formulation, derivative-informed training, and its use in evaluating the D-optimality rewards in \eqref{eq: incremental formulation immediate} and \eqref{eq: terminal formulation terminal}, culminating in a LANO-accelerated PGRL procedure for SBOED.

We first introduce dual dimension reduction for the parameter and state spaces to obtain low-dimensional representations suitable for learning. We then present the derivative-informed LANO architecture adjusted from \cite{go2025sequential} for random field parameter representing initial conditions, which predicts both state trajectories and sensitivities, describe the derivative-informed training objective, and finally explain how the trained surrogate accelerates Bayesian inversion and D-optimality reward evaluation within the RL loop.

\subsection{Dimension reduction}
\label{subsection: Dimension reduction}

As discussed in Section~\ref{subsection: High-fidelity discretization}, the random field parameter $m$ and the state variable $u_k(m)$ are discretized using the finite element method (FEM). These discretized representations, denoted by $\mathbf{m} \in \mathbb{R}^{d_m}$ and $\mathbf{u}_k \in \mathbb{R}^{d_u}$, are high-dimensional when a fine mesh is employed. LANO maps $\mathbf{m}$ and the initial condition to the state trajectory evaluated at a set of temporal points $\{t_{k'}\}_{k'=1}^{K'}$, which include the observation time grid $\{t_{k}\}_{k=1}^{K}$ on $[0,T]$. For notational convenience, we do not distinguish between these time grids. 

To mitigate the computational complexity associated with high-dimensional inputs and outputs in the neural operator, we employ dimension reduction techniques for both the parameter and state spaces.

\textbf{Active subspace (AS) projection for the parameter.} The dimensionality of the discretized parameter $\mathbf{m}$ is reduced via active subspace projection. In the context of PDE simulations, we use the expected cumulative Jacobian information of the state with respect to the parameter to define the following generalized eigenvalue problem:
\begin{align}
    \label{eq: active subspace}
    \mathbb{E}_{\mathbf{m}} \left[ \sum_{k=1}^{K} (\nabla_{\mathbf{m}} \mathbf{u}_k(\mathbf{m}))^\top \nabla_{\mathbf{m}} \mathbf{u}_k(\mathbf{m}) \right] \boldsymbol{\psi}^j_m = \lambda_j \mathbf{C}_{\text{prior}}^{-1} \boldsymbol{\psi}^j_m, \quad j = 1, \dots, r_m,
\end{align}
where $\lambda_1 \ge \dots \ge \lambda_{r_m}$ denote the $r_m$ largest generalized eigenvalues and $\{\boldsymbol{\psi}^j_m\}_{j=1}^{r_m}$ are the corresponding generalized eigenvectors, normalized such that $(\boldsymbol{\psi}^j_m)^\top \mathbf{C}_{\text{prior}}^{-1} \boldsymbol{\psi}^{j'}_m = \delta_{jj'}$. We approximate the expectation via Monte Carlo sampling from the prior \eqref{eq: prior}; the computation of the cumulative Jacobian information for each sample is discussed in Section~\ref{subsection: Jacobian computation for neural operator} for an example. The active subspace projection matrix is $\boldsymbol{\Psi}_m := [\boldsymbol{\psi}^1_m, \dots, \boldsymbol{\psi}^{r_m}_m] \in \mathbb{R}^{d_m \times r_m}$.
We define the dimensionally reduced parameter $\boldsymbol{\beta}_m \in \mathbb{R}^{r_m}$ via the projection:
\begin{align}
    \label{eq: parameter projection}
    \boldsymbol{\beta}_m := \boldsymbol{\Psi}^\top_m \mathbf{C}_{\text{prior}}^{-1} (\mathbf{m} - \mathbf{m}_{\text{prior}}).
\end{align}
Conversely, the parameter can be approximated as $\mathbf{m} \approx \mathbf{m}_{\text{prior}} + \boldsymbol{\Psi}_m \boldsymbol{\beta}_m$.

\textbf{PCA projection for the state.} The dimensionality of the state is reduced using principal component analysis (PCA). We generate $M_P$ samples of the random field parameter from the prior. For each sample $\mathbf{m}_i$, the time-dependent PDE \eqref{eq: pde} is solved via FEM to obtain the state trajectory $\{\mathbf{u}_k(\mathbf{m}^{(i)})\}_{k=1}^{K}$. These solutions are aggregated to form the snapshot matrix:
\begin{align}
    \label{eq: snapshot matrix}
    \mathbf{U} := [\mathbf{u}_{1}(\mathbf{m}^{(1)}), \dots, \mathbf{u}_K(\mathbf{m}^{(1)}), \dots, \mathbf{u}_{1}(\mathbf{m}^{(M_P)}), \dots, \mathbf{u}_K(\mathbf{m}^{(M_P)})] \in \mathbb{R}^{d_u \times K M_P}.
\end{align}
We center the snapshot matrix as $\hat{\mathbf{U}} = \mathbf{U} - \bar{\mathbf{u}}$, where $\bar{\mathbf{u}}$ is the empirical mean vector across all time steps and samples. Truncated Singular Value Decomposition (SVD) is applied to the centered matrix:
\begin{align}
    \label{eq: SVD}
    \hat{\mathbf{U}} \approx \hat{\mathbf{U}}_{r_u} := \boldsymbol{\Psi}_u \boldsymbol{\Sigma}_u \boldsymbol{\Phi}^\top_u,
\end{align}
where $\boldsymbol{\Sigma}_u = \text{diag}(\sigma_1, \dots, \sigma_{r_u})$ contains the $r_u$ largest singular values. The dimensionally reduced state variable $\boldsymbol{\beta}_u \in \mathbb{R}^{r_u}$ is obtained by projecting the centered state onto the principal basis $\boldsymbol{\Psi}_u$:
\begin{align}
    \label{eq: state projection}
    \boldsymbol{\beta}_u := \boldsymbol{\Psi}^\top_u (\mathbf{u} - \bar{\mathbf{u}}).
\end{align}
The reconstruction is given by $\mathbf{u} \approx \bar{\mathbf{u}} + \boldsymbol{\Psi}_u \boldsymbol{\beta}_u$.

\subsection{LANO architecture}
\label{subsection: LANO architecture}

The LANO architecture approximates the operator mapping the reduced parameter $\boldsymbol{\beta}_{m}$ and reduced initial condition $\boldsymbol{\beta}_{u_{0}}$ to the latent state sequence $\{ \boldsymbol{\beta}_{u_k} \}_{k=1}^{K}$. It consists of four stages: latent encoding, latent attention, latent dynamics, and decoding. Compared to \cite{go2025sequential} where the parameter represents the uncertain coefficient, we encode $\boldsymbol{\beta}_{m}$ and $\boldsymbol{\beta}_{u_{0}}$ to capture the parameter-dependent initial condition and use a transformer-style attention block with residual connections and layer normalization.

\textbf{Latent encoding.} We initialize a sequence of latent variables $\{ \mathbf{z}_{k} \}_{k=1}^{K}$ by combining a shared parameter projection with time-step specific state projections. Let $\mathbf{W}^{m}, \mathbf{b}^{m}$ be the global parameter projection weights, and $\mathbf{W}_{k}^{u_{0}}, \mathbf{b}_{k}^{u_{0}}$ be the time-specific initial condition weights. The encoding is defined as:
\begin{align}
\label{eq: latent encoding}
\mathbf{z}_{k} = \sigma_{z} \left( \mathbf{W}_{k}^{z} \left( \text{Concat}(\mathbf{W}^{m}\boldsymbol{\beta}_{m} + \mathbf{b}^{m}, \; \mathbf{W}_{k}^{u_{0}}\boldsymbol{\beta}_{u_{0}} + \mathbf{b}_{k}^{u_{0}}) \right) + \mathbf{b}_{k}^{z} \right) \in \mathbb{R}^{d_{z}},
\end{align}
for $k = 1, \dots, K$, where $\sigma_{z}$ is a non-linear activation. We denote the aggregate latent representation as $\mathbf{Z} := [ \mathbf{z}_{1}, \dots, \mathbf{z}_{K} ] \in \mathbb{R}^{d_{z} \times K}$. This use of time-specific projections allows the model to learn a temporal expansion of the parameter-dependent initial condition, effectively generating a trajectory of latent queries tailored to each simulation time step.

\textbf{Latent attention.} To capture temporal dependencies, we apply a causal self-attention mechanism. Positional embeddings are added to $\mathbf{Z}^\top$. We compute queries $\mathbf{Q}$, keys $\mathbf{K}$, and values $\mathbf{V}$ via linear projections of $\mathbf{Z}^\top$. The attention output is given by:
\begin{align}
\label{eq: attention score}
\mathcal{A} = \text{softmax} \left( \frac{\mathbf{Q}\mathbf{K}^{\top}}{\sqrt{d_a}} + \mathbf{M}_{\text{lt}} \right) \mathbf{V} \in \mathbb{R}^{K \times d_a},
\end{align}
where $\mathbf{M}_{\text{lt}}$ is a lower-triangular causal mask (with $-\infty$ in the upper triangle).

\textbf{Latent dynamics.} The attention output $\mathcal{A}$ is projected back to dimension $d_z$ and processed by a position-wise feed-forward network (FFN) with residual connections and layer normalization to produce the feature matrix $\mathbf{L} \in \mathbb{R}^{K \times d_z}$:
\begin{align}
\label{eq: transformed attention}
\mathbf{L} &= \text{LN}(\text{LN}(\mathbf{Z}^\top + \mathcal{A} \mathbf{W}^{\text{proj}}) + \text{FFN}(\text{LN}(\mathbf{Z}^\top + \mathcal{A} \mathbf{W}^{\text{proj}}))).
\end{align}
Subsequently, the reduced state variables are evolved autoregressively using a residual formulation. We simultaneously propagate the latent state $\boldsymbol{\beta}^{\mathbf{u}}_k$ and an auxiliary Jacobian latent representation $\boldsymbol{\beta}^{\mathbf{J}}_k$:
\begin{align}
\label{eq: latent dynamics state}
\boldsymbol{\beta}^{\mathbf{u}}_{k+1} &= \boldsymbol{\beta}^{\mathbf{u}}_{k} + \text{MLP}^{\mathbf{u}}_k(\boldsymbol{\ell}_{k+1}), \\
\label{eq: latent dynamics jacobian}
\boldsymbol{\beta}^{\mathbf{J}}_{k+1} &= \boldsymbol{\beta}^{\mathbf{J}}_{k} + \text{MLP}^{\mathbf{J}}_k(\boldsymbol{\ell}_{k+1}),
\end{align}
for $k = 0, \dots, K-1$, with initial conditions $\boldsymbol{\beta}^{\mathbf{u}}_{0} = \boldsymbol{\beta}^{\mathbf{J}}_{0} = \boldsymbol{\beta}_{u_{0}}$. Here, $\boldsymbol{\ell}_k$ is the $k$-th row of the feature matrix $\mathbf{L}$, and $\text{MLP}_k$ denotes a time-step specific multi-layer perceptron. The update increments depend on the attention features, while the autoregressive state enters through the residual addition.
The auxiliary sequence $\{\boldsymbol{\beta}^{\mathbf{J}}_k\}$ is introduced to enable accurate sensitivity prediction: its gradient with respect to the reduced parameter, $\nabla_{\boldsymbol{\beta}_m}\boldsymbol{\beta}^{\mathbf{J}}_{k}$, is used to reconstruct the state Jacobian in \eqref{eq: decoding jacobian}. It should be viewed as a learned latent representation for sensitivities rather than the Jacobian itself.

\textbf{Decoding.} The architecture bifurcates into two computational paths to reconstruct the high-fidelity quantities. The first path reconstructs the state $\mathbf{u}_k$ via the PCA basis $\boldsymbol{\Psi}_u$ defined in \eqref{eq: SVD}:
\begin{align}
\label{eq: decoding state}
\hat{\mathbf{u}}_k(\mathbf{m}) := \bar{\mathbf{u}} + \boldsymbol{\Psi}_u \boldsymbol{\beta}^{\mathbf{u}}_{k}.
\end{align}
The second path approximates the Jacobian of the state using the auxiliary variable $\boldsymbol{\beta}^{\mathbf{J}}_{k}$. The full Jacobian is reconstructed by projecting $\nabla_{\boldsymbol{\beta}_m}\boldsymbol{\beta}^{\mathbf{J}}_{k}$ using the bases $\boldsymbol{\Psi}_u$ and $\boldsymbol{\Psi}_m$:
\begin{align}
\label{eq: decoding jacobian}
\nabla_{\mathbf{m}}\hat{\mathbf{u}}_{k}(\mathbf{m}) := \boldsymbol{\Psi}_u \left(\nabla_{\boldsymbol{\beta}_m}\boldsymbol{\beta}^{\mathbf{J}}_{k} \right)\boldsymbol{\Psi}_m^\top \mathbf{C}_{\text{prior}}^{-1}.
\end{align}
The factor $\mathbf{C}_{\text{prior}}^{-1}$ arises from the whitening-based parameter projection in \eqref{eq: parameter projection}, which defines the reduced coordinates $\boldsymbol{\beta}_m$ in the prior-precision inner product.

\textbf{Implementation note.} In the training code, we use tanh in the latent encoding, fixed sinusoidal positional embeddings, a single attention head with $d_a = d_z/2$, and ELU activations in the per-step MLPs. The Jacobian term $\nabla_{\boldsymbol{\beta}_m}\boldsymbol{\beta}^{\mathbf{J}}_{k}$ is obtained by automatic differentiation of the Jacobian branch, and the reconstruction of $\nabla_{\mathbf{m}}\hat{\mathbf{u}}_{k}(\mathbf{m})$ in \eqref{eq: decoding jacobian} is applied outside the network and not used in the training.

\subsection{Data generation and derivative-informed training} 
\label{subsection: Data generation and derivative-informed training}

To construct the training dataset for the LANO surrogate, we sample $M_T$ realizations of the random field parameter from the prior distribution \eqref{eq: prior}, denoted as $\{ \mathbf{m}^{(i)} \}_{i=1}^{M_T}$. For each realization $\mathbf{m}^{(i)}$, the time-dependent PDE \eqref{eq: pde} is solved via the finite element method to obtain the high-fidelity state trajectory $\{ \mathbf{u}_k(\mathbf{m}^{(i)}) \}_{k=1}^{K}$ and the associated Jacobians. Computing these Jacobians is nontrivial; we detail the procedure in Section~\ref{subsection: Jacobian computation for neural operator} for an example. These high-dimensional quantities are projected onto the corresponding reduced bases to form the training targets. Let $\boldsymbol{\beta}_m^{(i)}$ denote the projected parameter. The projected state $\boldsymbol{\beta}_{\mathbf{u}, k}^{(i)}$ and projected Jacobian $\boldsymbol{\beta}_{\mathbf{J}, k}^{(i)}$ are computed as:
\begin{align}
    \label{eq: data projection}
    \boldsymbol{\beta}_{\mathbf{u}, k}^{(i)} = \boldsymbol{\Psi}_u^\top (\mathbf{u}_k(\mathbf{m}^{(i)}) - \bar{\mathbf{u}}), \quad \boldsymbol{\beta}_{\mathbf{J}, k}^{(i)} = \boldsymbol{\Psi}_u^\top \nabla_{\mathbf{m}} \mathbf{u}_k(\mathbf{m}^{(i)}) \boldsymbol{\Psi}_m.
\end{align}

We formulate a derivative-informed loss function:
\begin{align}
    \label{eq: training loss}
    \mathcal{L}_{\text{LANO}}(\boldsymbol{\Theta}) := \frac{1}{M_T} \sum_{i=1}^{M_T} \sum_{k=1}^{K} \left( \| \boldsymbol{\beta}_{\mathbf{u}, k}^{(i)} - \hat{\boldsymbol{\beta}}^{\mathbf{u}}_{k}(\boldsymbol{\beta}_m^{(i)}; \boldsymbol{\Theta}) \|_2^2 + \lambda \| \boldsymbol{\beta}_{\mathbf{J}, k}^{(i)} - \nabla_{\boldsymbol{\beta}_m} \hat{\boldsymbol{\beta}}^{\mathbf{J}}_{k}(\boldsymbol{\beta}_m^{(i)}; \boldsymbol{\Theta}) \|_F^2 \right),
\end{align}
where $\boldsymbol{\Theta}$ represents the trainable parameters of the LANO architecture. This loss minimizes the discrepancy between the network predictions and the projected targets for both the state and Jacobian, thereby enforcing consistency in the solution and its sensitivity \cite{olearyroseberry2024derivativeinformed,qiu2024derivativeenhanced,go2025sequential}. Here, $\hat{\boldsymbol{\beta}}^{\mathbf{u}}_{k}$ and $\hat{\boldsymbol{\beta}}^{\mathbf{J}}_{k}$ denote the outputs of the state computational path \eqref{eq: latent dynamics state} and the Jacobian computational path \eqref{eq: latent dynamics jacobian}, respectively. $\nabla_{\boldsymbol{\beta}_m} \hat{\boldsymbol{\beta}}^{\mathbf{J}}_{k}$ is computed by automatic differentiation of the Jacobian branch. The hyperparameter $\lambda > 0$ balances the two terms.

\subsection{Acceleration of solving Bayesian inverse problems} \label{subsection: Acceleration of solving Bayesian inverse problems}

The computational cost of the SBOED framework is dominated by repeatedly solving the Bayesian inverse problem. We accelerate this process by replacing the high-fidelity FEM solver with the trained LANO surrogate within the Laplace approximation framework described in Section~\ref{subsection: Laplace approximation of the posterior distribution}.

For the MAP estimation, the PtO maps in the optimization problem \eqref{eq: MAP optimization} are approximated by the LANO predictions. The optimization is performed with respect to $\boldsymbol{\beta}_m$ in the reduced parameter space:
\begin{align}
    \label{eq: accelerated MAP}
    \min_{\boldsymbol{\beta}_m \in \mathbb{R}^{r_m}} \frac{1}{2} \sum_{\ell=1}^{n} \sum_{k \in \mathcal{K}_\ell} \| \mathbf{y}^{(k)} - \mathbf{B}_\ell(\bar{\mathbf{u}} + \boldsymbol{\Psi}_u \boldsymbol{\beta}^{\mathbf{u}}_{k}(\boldsymbol{\beta}_m)) \|^2_{\boldsymbol{\Gamma}_{\text{noise}}^{-1}} + \frac{1}{2} \| \boldsymbol{\beta}_m \|^2.
\end{align}
Note that $\mathbf{B}_\ell(\bar{\mathbf{u}} + \boldsymbol{\Psi}_u \boldsymbol{\beta}^{\mathbf{u}}_{k}(\boldsymbol{\beta}_m))$ can be efficiently evaluated for linear observation operators (pointwise sensors) $\mathbf{B}_\ell \in \mathbb{R}^{d_y \times d_u}$ as $\mathbf{B}_\ell(\bar{\mathbf{u}} + \boldsymbol{\Psi}_u \boldsymbol{\beta}^{\mathbf{u}}_{k}(\boldsymbol{\beta}_m)) = \mathbf{B}_\ell \bar{\mathbf{u}} + (\mathbf{B}_\ell \boldsymbol{\Psi}_u) \boldsymbol{\beta}^{\mathbf{u}}_{k}(\boldsymbol{\beta}_m)$ where $\mathbf{B}_\ell \bar{\mathbf{u}} \in \mathbb{R}^{d_y}$ and $\mathbf{B}_\ell \boldsymbol{\Psi}_u \in \mathbb{R}^{d_y \times r_u}$ can be pre-computed once. The regularization term simplifies to the Euclidean norm due to the whitening transformation in \eqref{eq: parameter projection}. Let $\boldsymbol{\beta}_{\text{MAP}}^{\mathbf{x}_{n+1}}$ denote the solution, latent MAP, to this reduced problem.

For the computation of the posterior covariance, the Jacobian of the misfit term is required. We utilize LANO's second computational path \eqref{eq: decoding jacobian} to approximate the Gauss-Newton Hessian \eqref{eq: GN approximation}:
\begin{align}
    \label{eq: accelerated Hessian}
    \mathbf{H}_{\text{GN}}^{\mathbf{x}_{n+1}} & \approx \sum_{\ell=1}^{n} \sum_{k \in \mathcal{K}_\ell} \left( \nabla_{\mathbf{m}} \mathbf{B}_\ell(\hat{\mathbf{u}}_{k})\right)^\top \boldsymbol{\Gamma}_{\text{noise}}^{-1} \left( \nabla_{\mathbf{m}} \mathbf{B}_\ell(\hat{\mathbf{u}}_{k}) \right),
\end{align}
where $\nabla_{\mathbf{m}} \mathbf{B}_\ell(\hat{\mathbf{u}}_{k}) = \mathbf{B}_\ell \boldsymbol{\Psi}_u \left(\nabla_{\boldsymbol{\beta}_m} \boldsymbol{\beta}_k^{\mathbf{J}} \right)\boldsymbol{\Psi}_m^\top \mathbf{C}_{\text{prior}}^{-1}$ for linear observation operator $\mathbf{B}_\ell$. This enables efficient pre-computation of the small matrix $\mathbf{R}_{\ell} = \boldsymbol{\Psi}_u^\top \mathbf{B}_\ell^\top \boldsymbol{\Gamma}_{\text{noise}}^{-1} \mathbf{B}_\ell \boldsymbol{\Psi}_u \in \mathbb{R}^{r_u\times r_u}$.
Then the generalized eigenvalue problem \eqref{eq: generalized eigenvalue} can be transformed as the eigenvalue problem in the reduced parameter space:
\begin{align}\label{eq: reduced geneig}
    \mathbf{\hat H}_{\text{GN}}^{\mathbf{x}_{n+1}} (\boldsymbol{\beta}_{\text{MAP}}^{\mathbf{x}_{n+1}}) \mathbf{v}_j = \lambda_j \mathbf{v}_j, \quad j = 1, \dots, r_m,
\end{align}
where the eigenvectors of \eqref{eq: generalized eigenvalue} can be approximated as $\mathbf{w}_j = \boldsymbol{\Psi}_m \mathbf{v}_j$ and the reduced Hessian $\mathbf{\hat H}_{\text{GN}}^{\mathbf{x}_{n+1}} (\boldsymbol{\beta}_{\text{MAP}}^{\mathbf{x}_{n+1}}) \in \mathbb{R}^{r_m\times r_m}$ can be efficiently computed as  
\begin{align}
    \mathbf{\hat H}_{\text{GN}}^{\mathbf{x}_{n+1}} (\boldsymbol{\beta}_{\text{MAP}}^{\mathbf{x}_{n+1}}) = \sum_{\ell=1}^{n} \sum_{k \in \mathcal{K}_\ell} \left(\nabla_{\boldsymbol{\beta}_m} \boldsymbol{\beta}_k^{\mathbf{J}}(\boldsymbol{\beta}_{\text{MAP}}^{\mathbf{x}_{n+1}})\right)^\top \mathbf{R}_{\ell} \nabla_{\boldsymbol{\beta}_m} \boldsymbol{\beta}_k^{\mathbf{J}}(\boldsymbol{\beta}_{\text{MAP}}^{\mathbf{x}_{n+1}}).
\end{align}
This leads directly to the D-optimality computation in \eqref{eq: D-optimality computation} using the $r_m$ eigenvalues from \eqref{eq: reduced geneig}. 

To further reduce the computational cost, we can replace the reduced MAP estimate $\boldsymbol{\beta}_{\text{MAP}}^{\mathbf{x}_{n+1}}$ with the prior sample that is used in generating the corresponding observation data for the MAP estimate, which avoids solving the optimization problem \eqref{eq: accelerated MAP} altogether. This approximation is reasonable when the D-optimality criteria evaluated at the prior sample and the MAP point are highly correlated, as exploited in \cite{wu2023fast} in the full space and observed in our numerical experiments.

\subsection{Pseudocode for the overall algorithm}

We summarize the components above in Algorithm~\ref{alg: pgrl}, which applies PGRL (Section~\ref{section: Policy gradient reinforcement learning}) to SBOED with infinite-dimensional random field parameters. The scalable approximations in Section~\ref{section: Scalable approximations for SBOED} are further accelerated using the LANO surrogate developed in this section. We note that Steps~\ref{step: sample} and~\ref{step: surrogate}, which do not depend on the experimental designs, can be precomputed for all Monte Carlo samples of the prior to further reduce wall-clock time. Moreover, the projected prior samples $\boldsymbol{\beta}_m^{(i)}$ in Step~\ref{step: sample} can be drawn directly from a standard normal distribution due to the whitening transformation in \eqref{eq: parameter projection}.
\begin{algorithm}[!htb]
\caption{LANO-accelerated PGRL for SBOED}
\label{alg: pgrl}
\begin{algorithmic}[1]
\STATE Initialize algorithmic hyperparameters: number of policy updates $L$, Q-function updates $L'$, Monte Carlo sample size $M$. Define architectures for policy network $\mathcal{P}_{\mathbf{w}}$ and Q-network $\mathcal{Q}_{\boldsymbol{\eta}}$;
\STATE Initialize network parameters $\mathbf{w}$ and $\boldsymbol{\eta}$.
\FOR{$l = 1, \cdots, L$}
    \STATE \textbf{Data Collection (Surrogate Simulation):}
    \FOR{$i = 1, \cdots, M$}
        \STATE \label{step: sample} Sample random field parameter $\mathbf{m}^{(i)}$ from the prior \eqref{eq: prior} and compute its projection $\boldsymbol{\beta}_m^{(i)}$ via \eqref{eq: parameter projection}; 
        \STATE \label{step: surrogate} Compute state trajectory approximations $\{\hat{\mathbf{u}}_k(\mathbf{m}^{(i)})\}_{k=1}^K$ via the trained LANO surrogate \eqref{eq: decoding state}; 
        \STATE Initialize state history $\mathbf{x}_1^{(i)} = \emptyset$ and its zero-padded representation \eqref{eq: state variable neural network input};
        \FOR{$n = 1, \cdots, N$}
            \STATE Compute design vector $\mathbf{d}_n = \mathcal{P}_{\mathbf{w}}(n, \mathbf{x}_n^{(i)})$ and its associated observation operator $\mathbf{B}_n$;
            \STATE Generate synthetic observations $\mathbf{y}_n$ by applying $\mathbf{B}_n$ to the surrogate states $\{\hat{\mathbf{u}}_k(\mathbf{m}^{(i)})\}_{k=1}^K$;
            \STATE Update the state $\mathbf{x}_{n+1}^{(i)} \leftarrow [\mathbf{x}_n^{(i)}, \mathbf{d}_n, \mathbf{y}_n]$ and its zero-padded representation \eqref{eq: state variable neural network input};
        \ENDFOR
\STATE Compute terminal rewards $\{g_n\}_{n=1}^{N+1}$ \eqref{eq: terminal formulation immediate} and \eqref{eq: terminal formulation terminal}) using the D-optimality criterion \eqref{eq: D-optimality computation}, where the eigenvalues are computed as in \eqref{eq: reduced geneig} accelerated by LANO (Section~\ref{subsection: Acceleration of solving Bayesian inverse problems}).
    \ENDFOR
    \STATE \textbf{Critic Update:}
    \FOR{$l' = 1, \cdots, L'$}
        \STATE Update Q-network parameters $\boldsymbol{\eta}$ by minimizing \eqref{eq: Q-network training loss function} via stochastic gradient descent. 
    \ENDFOR
    \STATE \textbf{Actor Update:}
    \STATE Update policy parameters $\mathbf{w}$ via stochastic gradient ascent using the estimator \eqref{eq: MC policy gradient}.
\ENDFOR
\STATE \textbf{Return} optimized policy $\pi_{\mathbf{w}}$.
\end{algorithmic}
\end{algorithm}


\section{Case study: Optimal tracking of contaminant source} 
\label{section: Case study: contaminant source tracking}

\textbf{Case-study roadmap.} The unknown parameter is a random field $m(x)$ that enters the model of contaminant transport through the initial condition $u(0,x)=m^2(x)$. We observe the evolving concentration field $u(t,x)$ through sparse pointwise sensors; the design variables are the sensor locations, encoded by stage-wise displacement vectors $\{\mathbf{d}_n\}_{n=1}^N$ that induce linear observation operators $\{\mathbf{B}_n\}_{n=1}^N$. We perform Bayesian inversion for $m$ from noisy observations and evaluate information gain using a terminal D-optimality reward under a Laplace approximation. We then apply the learned policy to adapt sensor placements sequentially.

\subsection{Time-dependent advection-diffusion equation with random initial condition}

\begin{figure}[htbp]
  \centering
  
  \begin{subfigure}[t]{0.31\linewidth}
    \centering
    \includegraphics[width=\linewidth]{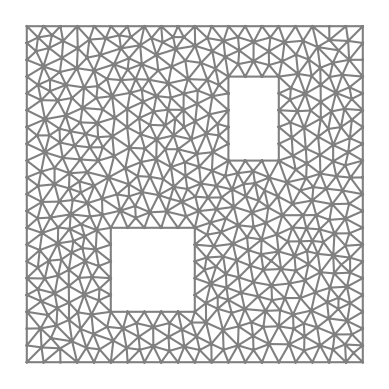}
    \caption{$\Omega$ and mesh}
    \label{fig: physical domain + velocity fields a}
  \end{subfigure}
  \hspace{0.25em}
  \begin{subfigure}[t]{0.31\linewidth}
    \centering
    \includegraphics[width=\linewidth]{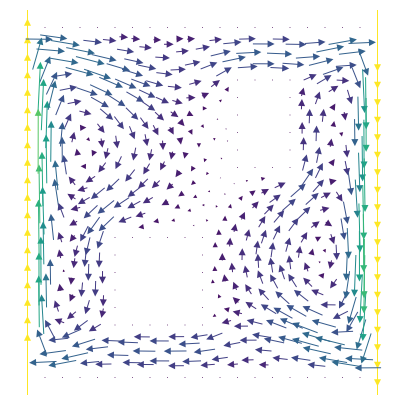}
    \caption{$\vec{v}$ with $\vec{g}_{1}$}
    \label{fig: physical domain + velocity fields b}
  \end{subfigure}
  \hspace{0.25em}
  \begin{subfigure}[t]{0.31\linewidth}
    \centering
    \includegraphics[width=\linewidth]{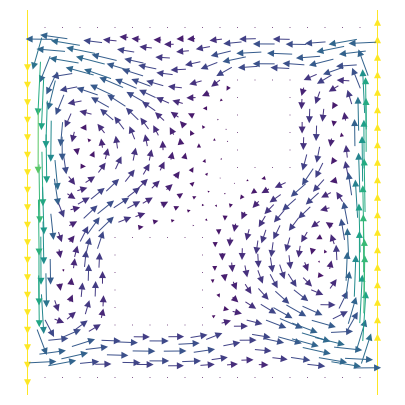}
    \caption{$\vec{v}$ with $\vec{g}_{2}$}
    \label{fig: physical domain + velocity fields c}
  \end{subfigure}
    
  \caption{(a) The physical domain $\Omega \subset \mathbb{R}^{2}$ featuring two impenetrable buildings (blocks) and the computational mesh (refined twice in our computation). (b) The steady-state velocity field $\vec{v}$ resulting from boundary condition $\vec{g}_{1}$. (c) The velocity field resulting from boundary condition $\vec{g}_{2}$, which induces a reversed flow pattern.}
  \label{fig: physical domain + velocity fields}
\end{figure}

We consider a variant of the time-dependent advection-diffusion equation from \cite{petra2011model, villa2021hippylib}, which models the diffusion and transport of a contaminant driven by a fluid flow. We define the initial condition as the square of a Gaussian random field parameter, $m^2$, where the physical domain $\Omega \subset \mathbb{R}^{2}$ includes two buildings that are impenetrable to the contaminant (Figure~\ref{fig: physical domain + velocity fields a}). Squaring the parameter $m$ ensures that the initial contaminant concentration remains non-negative. The governing equation is:
\begin{align}
    \frac{\partial u}{\partial t} - \kappa\Delta u + \vec{v}\cdot \nabla u = 0, & \text{ in } \Omega \times (0, T), \nonumber \\
    \label{eq: time-dependent advection-diffusion equation}
    u(0, \cdot) = m^{2}(\cdot), & \text{ in } \Omega, \\
    \nabla u \cdot \vec{n} = 0, & \text{ on } \partial\Omega \times (0, T), \nonumber
\end{align}
where $T=4.0$ is the terminal time, $\kappa =0.001$ is the diffusivity coefficient, and the state variable $u(t, x)$ represents the contaminant concentration at time $t$ and location $x$. Our goal is twofold: (i) infer the unknown random field parameter $m$ from sparse pointwise sensor observations via Bayesian inversion (Section~\ref{subsection: Bayesian inverse problem given sparse sensor observations}), and (ii) design sensor placements to maximize expected information gain via SBOED (Section~\ref{subsection: Optimal sensor placement}).

In \eqref{eq: time-dependent advection-diffusion equation}, the velocity field $\vec{v}$ is determined by the steady-state Navier--Stokes equations:
\begin{align}
    -\frac{1}{\text{Re}}\Delta\vec{v} + \nabla q + \vec{v}\cdot\nabla\vec{v} = 0, & \text{ in } \Omega, \nonumber \\
    \label{eq: steady-state Navier-Stokes}
    \nabla \cdot \vec{v} = 0, & \text{ in } \Omega, \\
    \vec{v} = \vec{g}, & \text{ on } \partial\Omega, \nonumber
\end{align}
where $q$ denotes pressure and $\text{Re}=100$ is the Reynolds number. We consider two sets of Dirichlet boundary conditions, $\vec{g}_{1}$ and $\vec{g}_{2}$. For $\vec{g}_{1}$, we set $\vec{v} = \vec{e}_{2}$ on the left wall, $\vec{v} = -\vec{e}_{2}$ on the right wall, and $\vec{v} = \vec{0}$ on the top and bottom walls, where $\vec{e}_{2}$ is the unit vector along the vertical axis. The condition $\vec{g}_{2}$ is obtained by negating $\vec{g}_{1}$ on the left and right walls. The resulting velocity fields are visualized in Figures~\ref{fig: physical domain + velocity fields b} and~\ref{fig: physical domain + velocity fields c}. Unless otherwise specified, all numerical experiments use the velocity field determined by $\vec{g}_{1}$.

\subsection{Bayesian inversion with sparse sensor observations} \label{subsection: Bayesian inverse problem given sparse sensor observations}

\begin{figure}[!htb]
  \centering  
    \includegraphics[width=\linewidth]{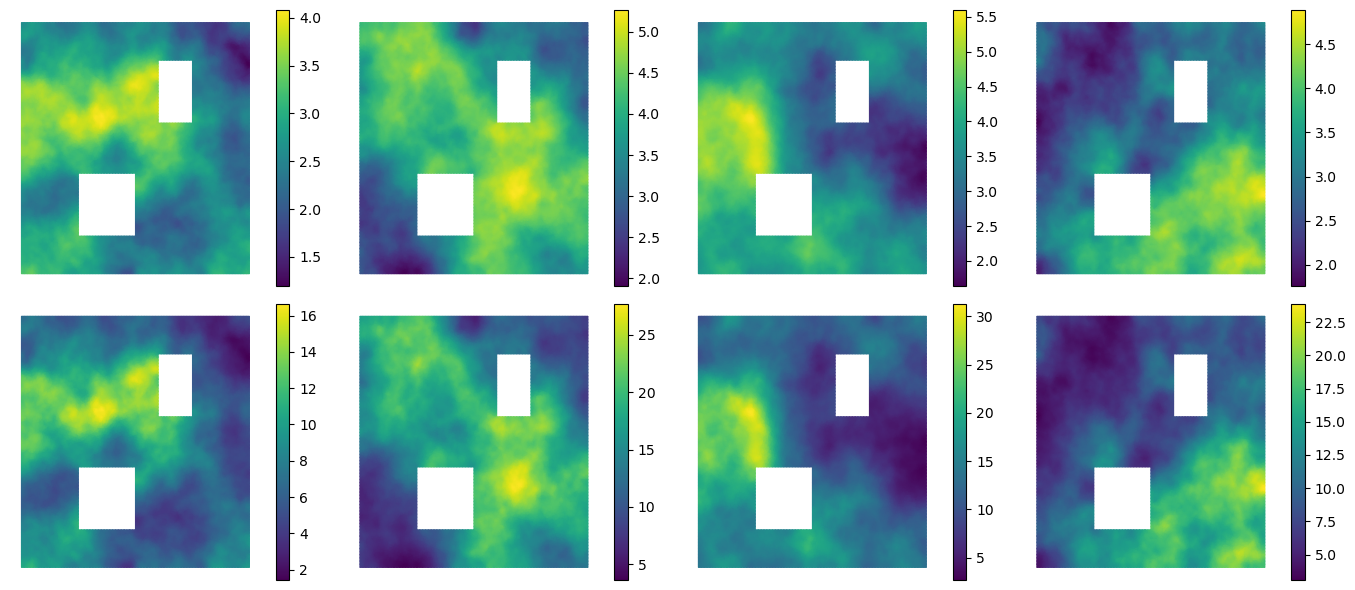}
  \caption{Top: Four independent samples drawn from the Gaussian prior distribution of the parameter $m$. The selected covariance parameters ($\gamma, \delta$) induce significant spatial variability, resulting in distinct high-value regions across samples. Bottom: The corresponding initial condition fields, obtained by squaring the prior samples ($m^2$). }
  \label{fig: prior samples + initial condition samples}
\end{figure}

The prior for the random field parameter $m$ in the initial condition of \eqref{eq: time-dependent advection-diffusion equation} is a Gaussian random field \eqref{eq: prior} with Robin boundary condition $\gamma \nabla m \cdot \vec{n} + \beta m = 0$ on $\partial \Omega$. Following \cite{villa2021hippylib, daon2018mitigating, roininen2014whittlematern}, the Robin coefficient is set to $\beta=\frac{\gamma \sqrt{\delta / \gamma}}{1.42}$ to mitigate boundary artifacts. The mean field is constant with $m_{\text{prior}} = 3.0$. We use $\alpha =2$, $\gamma =0.1$, and $\delta =0.8$, so that $\mathcal{C}_{\text{prior}} = \mathcal{A}^{-2} = (-0.1 \Delta + 0.8I)^{-2}$. These values induce significant spatial variability in the prior samples and the resulting initial condition $m^2$; see Figure~\ref{fig: prior samples + initial condition samples}.

We consider sparse sensor observations collected over $K=19$ time steps, starting at $t=0.4$ and ending at $T=4.0$ with a time interval of $0.2$. We conduct $N=4$ experiments, placing sensors at distinct locations in each experiment to observe the state variable. The observation time indices are grouped as follows:
\begin{align}\label{eq: sets of indices of observation time points}
    \mathcal{K}_{1} = \{ 1, 2 \}, \quad \mathcal{K}_{2} = \{ 3, 4 \}, \quad \mathcal{K}_{3} = \{ 5, 6 \}, \quad \mathcal{K}_{4} = \{ 7, \cdots, 19 \}.
\end{align}
Equivalently, the observation times are $t_k = 0.4 + 0.2(k-1)$ for $k=1,\dots,19$, and sensor placements are updated only at the stage boundaries (i.e., at $t=0.4, 0.8, 1.2,$ and $1.6$), remaining fixed for all $k \in \mathcal{K}_n$ within each experiment $n$.

\begin{figure}[!htb]
  \centering
  \includegraphics[width=\linewidth]{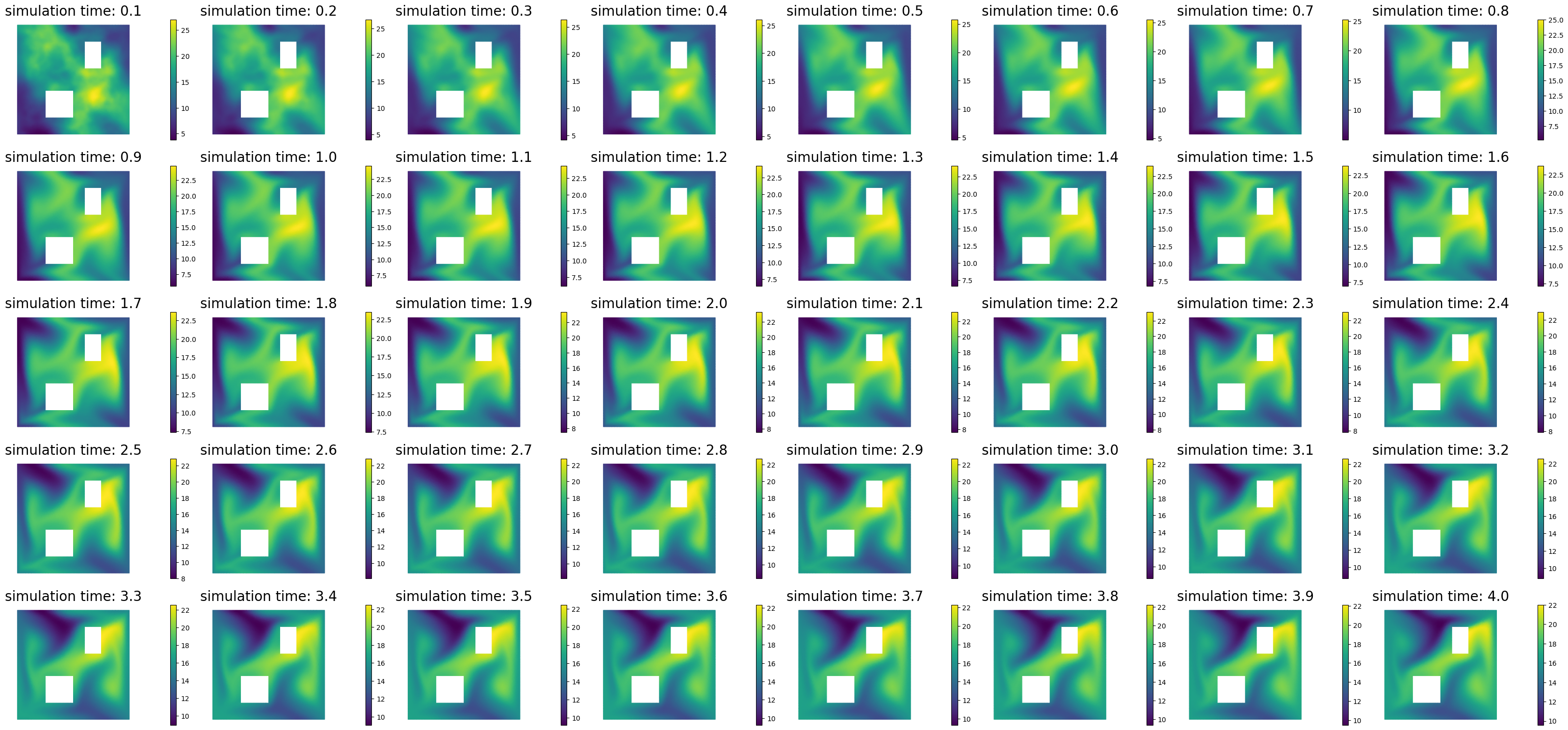}
  \caption{Time evolution of the contaminant concentration $u(t, x)$, simulated using the Finite Element Method (FEM) for a parameter $m$ sampled from the prior. The state variable exhibits negligible variation after $t=1.6$, which marks the final sensor placement time step. This steady-state behavior justifies the concentration of experiments in the earlier time intervals.}
  \label{fig: FEM-trajectory}
\end{figure}

To simulate \eqref{eq: time-dependent advection-diffusion equation}, we employ $K'=41$ simulation time steps spanning from $t=0.0$ to $T=4.0$ with a step size of $0.1$, i.e., $t_{k'} = 0.1(k'-1)$ for $k' = 1,\dots,41$. We concentrate experiments during the initial observation windows because the temporal variation of the state variable diminishes over time; see Figure~\ref{fig: FEM-trajectory}.

The dimension of the observation vectors matches the number of sensors when utilizing a pointwise observation model. After discretization, the observation operators $\{\mathbf{B}_{n}\}_{n=1}^N$ are matrices of size $d_y \times d_u$. The observation noise covariance matrix $\boldsymbol{\Gamma}_{\text{noise}}$ for $\{\boldsymbol{\epsilon}_{k}\}_{k=1}^K$ is diagonal, with entries set to the square of 1\% of a reference observation magnitude $y_{\max} := \max_{k,s} |y^{(k)}_s|$ (maximum over all sensors and observation times in a trajectory). We denote this diagonal entry by $\sigma^2$.

\subsection{Setup for optimal sensor placement} \label{subsection: Optimal sensor placement}

\begin{figure}[!htb]
    \centering
    \begin{subfigure}[t]{0.31\linewidth}
        \centering
        \includegraphics[width=\linewidth]{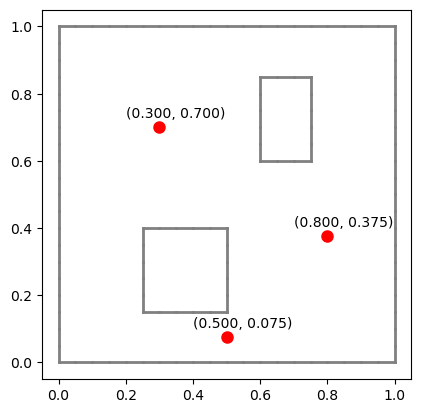}
        \caption{Initial sensor placement}
        \label{fig: initial placement of sensors}
    \end{subfigure}
    \hfill
    \begin{subfigure}[t]{0.31\linewidth}
        \centering
        \includegraphics[width=\linewidth]{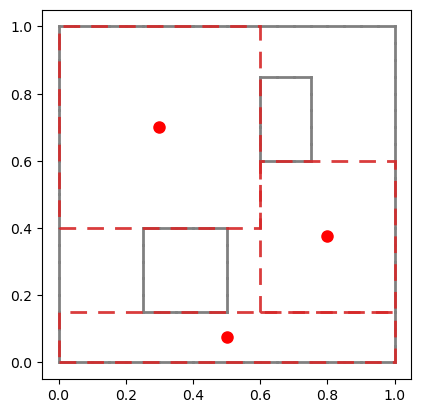}
        \caption{Permissible sensor regions}
        \label{fig: regions sensors can reside}
    \end{subfigure}
    \hfill
    \begin{subfigure}[t]{0.31\linewidth}
        \centering
        \includegraphics[width=\linewidth]{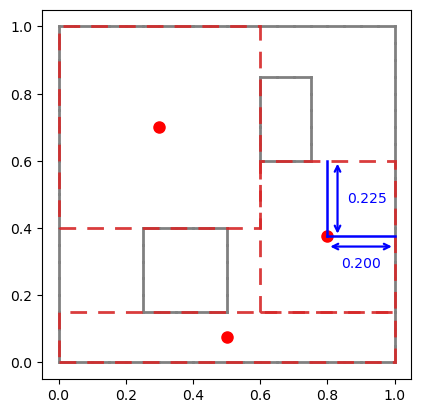}
        \caption{Distance to boundaries}
        \label{fig: regions sensors can reside distance}
    \end{subfigure}
    \caption{(a) The initial positions of the three sensors used for sparse observations. (b) The red dashed rectangles indicate the permissible regions within which each sensor must remain. (c) An illustration of the distance constraints for the sensor at $(0.800, 0.375)$; the bounds for the design vector are derived from these distances to ensure validity across all $N$ experiments.}
    \label{fig: sensor placement}
\end{figure}

We employ a set of three sensors to collect sparse observations, with their initial configurations shown in Figure~\ref{fig: initial placement of sensors}. Let $s \in \{1, 2, 3\}$ index the sensors, and let $(x_1[s], x_2[s])$ denote the coordinates of the $s$-th sensor. In each experiment $n$, the design vector $\mathbf{d}_n \in \mathbb{R}^{6}$ defines the displacement applied to the current sensor locations. The update rule for the $s$-th sensor in the $n$-th experiment is:
\begin{align}
    \label{eq: sensor placement update x coordinate}
    & x_1[s] \leftarrow x_1[s] + \mathbf{d}_{n}[2s - 1], \\
    \label{eq: sensor placement update y coordinate}
    & x_2[s] \leftarrow x_2[s] + \mathbf{d}_{n}[2s],
\end{align}
where $\mathbf{d}_{n}[k]$ represents the $k$-th entry of the design vector $\mathbf{d}_{n}$. 

To ensure that the sensors remain within the physical domain $\Omega$ throughout the sequential experiments, we impose constraints on the design vector. We define a permissible rectangular region for each sensor (Figure~\ref{fig: regions sensors can reside}). The step size limits for the design vector are derived by dividing the distance from the initial coordinates to the boundaries of these regions by the total number of experiments, $N$. For example, for the sensor initially located at $(0.800, 0.375)$, the distances to the nearest boundaries are $0.200$ and $0.225$ (Figure~\ref{fig: regions sensors can reside distance}). Therefore, to guarantee the sensor stays within bounds after $N=4$ experiments, the magnitude of the corresponding entries in $\mathbf{d}_n$ is constrained to be at most $0.200/N$ and $0.225/N$ for simplicity.

\subsection{Computation of the gradient and Hessian of the cost functional}\label{sec:gradHess}
To facilitate the inexact Newton-CG algorithm for computing the MAP point by solving \eqref{eq: MAP optimization}, we require the gradient and Hessian of the cost function. Because the squared initial condition renders the PtO map nonlinear with respect to the random field parameter, we provide the derivations below. Readers primarily interested in the overall algorithmic pipeline may skim this subsection; the derivations here explain how the gradient and Hessian actions are obtained for MAP computation and low-rank eigenvalue evaluations.

In the continuous setting (prior to discretization), the cost functional $\mathcal{J}(m)$ from \eqref{eq: MAP optimization} is expressed as:
\begin{align}
    \label{eq: MAP optimization cost functional integral}
    \mathcal{J}(m) =& \frac{1}{2\sigma^{2}} \sum_{\ell=1}^{n}\sum_{k \in \mathcal{K}_{\ell}} \int_{0}^{T} \left( \mathcal{B}_\ell u - \mathbf{y}^{(k)} \right)^{T} \left( \mathcal{B}_\ell u - \mathbf{y}^{(k)} \right) \delta(t - t_k) dt + \frac{1}{2} \int_{\Omega} \left| \mathcal{A} \left( m - m_{\text{prior}} \right)(x) \right|^{2} dx,
\end{align}
where $\mathcal{B}_\ell$ is the linear observation operator for the $\ell$-th experiment, active at time instances $t_k$ for $k \in \mathcal{K}_\ell$. The operator $\mathcal{A}$ is the elliptic differential operator defining the prior covariance $\mathcal{C}_{\text{prior}} = \mathcal{A}^{-2}$ (with $\alpha=2$).

Let $L^2(\Omega)$ denote the space of square-integrable functions over $\Omega$, and $H^1(\Omega)$ denote the space of functions in $L^2(\Omega)$ whose first derivatives are also square-integrable. The state variable satisfies $u(t, \cdot) \in H^1(\Omega)$ for $t \in (0, T)$. Accordingly, we define the trial space $\mathcal{U}$ and test space $\mathcal{P}$ as:
\begin{align}
    \label{eq: trial function space}
    & \mathcal{U} := L^{2}\left( \left( 0, T \right); H^{1}(\Omega) \right) \cap C^{0}\left( \left( 0, T \right); L^{2}(\Omega) \right), \\
    \label{eq: test function space}
    & \mathcal{P} := L^{2}\left( \left( 0, T \right); H^{1}(\Omega) \right) \cap C^{0}\left( \left( 0, T \right); L^{2}(\Omega) \right).
\end{align}
To derive the weak form of the state equation \eqref{eq: time-dependent advection-diffusion equation}, we multiply the residual by a test function $p\in \mathcal{P}$ and integrate over the space-time domain. We enforce the initial condition weakly using a test function $q(x) \in H^1(\Omega)$. Adding these terms to the cost functional \eqref{eq: MAP optimization cost functional integral} yields the Lagrangian functional $\mathcal{L}^g$:
\begin{align}\label{eq: Lagrangian functional gradient}
    \mathcal{L}^{g}(u, p, q, m) :=& \mathcal{J}(m) + \int_{0}^{T}\int_{\Omega} \left( \frac{\partial u}{\partial t} p + \kappa \nabla u \cdot \nabla p + (\vec{v} \cdot \nabla u) p \right) dxdt + \int_{\Omega}q(x) \left( u(0, x) - m^{2}(x) \right) dx,
\end{align}
where $(u, p, q, m) \in \mathcal{U} \times \mathcal{P} \times H^{1}(\Omega) \times \mathcal{M}$, with $\mathcal{M} = H^{1}(\Omega)$ denoting the parameter space.

{\bf Compute the gradient.}
We compute the gradient of $\mathcal{J}(m)$ with respect to $m$ using the adjoint method. First, given $m$, we solve the forward problem \eqref{eq: time-dependent advection-diffusion equation} to obtain the state $u \in \mathcal{U}$. Next, we require the first variation of $\mathcal{L}^g$ with respect to $u$ to vanish for all perturbations $\hat{u} \in \mathcal{U}$:
\begin{align}\label{eq: adjoint problem weak form}
    0 = \delta_{u}\mathcal{L}^{g} \hat{u} 
    = & \frac{1}{\sigma^{2}}\sum_{\ell=1}^{n}\sum_{k \in \mathcal{K}_{\ell}}\int_{0}^{T} \left( \mathcal{B}_\ell u - \mathbf{y}^{(k)} \right)^{T} \left( \mathcal{B}_\ell\hat{u} \right) \delta(t - t_k) dt \nonumber \\
    & + \int_{0}^{T}\int_{\Omega} \left( - \hat{u} \frac{\partial p}{\partial t} + \kappa \nabla \hat{u} \cdot \nabla p + (\vec{v} \cdot \nabla \hat{u}) p \right) dxdt,
\end{align}
where we used integration by part and set the terminal condition $p(T, x)=0$ and $q(x) = p(0,x)$.
This equation constitutes the weak form of the adjoint problem, which determines the adjoint variable $p$ by solving the adjoint problem backward in time. 

The weak form of the gradient, denoted by $(G(m), \hat{m})$, is derived from the variation of $\mathcal{L}^g$ with respect to $m$ in an arbitrary direction $\hat{m} \in \mathcal{M}$:
\begin{align}\label{eq: cost functional gradient evaluation in direction}
    &\left( G(m), \hat{m} \right) = \delta_{m}\mathcal{L}^{g} \hat{m}  
    = \int_{\Omega} \mathcal{A}\left( m - m_{\text{prior}} \right)(x) \mathcal{A}(\hat{m})(x) dx - \int_{\Omega}2 p(0, x)m(x)\hat{m}(x) dx.
\end{align}
Eliminating the test function $\hat{m}$, we obtain the strong form of the gradient $G(m)$:
\begin{equation}
\label{eq: cost functional gradient evaluation}
    G(m)=
        \begin{cases}
        \mathcal{A}^{2}(m - m_{\text{prior}}) - 2p(0,\cdot)m,
        & \text{in } \; \Omega, \\
        \gamma\nabla m \cdot \vec{n} + \beta m,
        & \text{on } \; \partial\Omega.
        \end{cases}
\end{equation}

{\bf Compute the Hessian action.}
Next, we evaluate the action of the Hessian of $\mathcal{J}(m)$ in a direction $\tilde{m}$. We define the second-order Lagrangian functional $\mathcal{L}^{h}$. The first term corresponds to the weak form of the gradient \eqref{eq: cost functional gradient evaluation in direction}, with the test function $\hat{m}$ replaced by the direction $\tilde{m}$. The subsequent terms represent the weak forms of the state equation \eqref{eq: time-dependent advection-diffusion equation} (with the adjoint variable $p$ replaced by the incremental adjoint $\tilde{p}$) and the adjoint problem (with the state variation $\hat{u}$ replaced by the incremental state $\tilde{u}$).
\begin{align}\label{eq: Lagrangian functional Hessian}
\mathcal{L}^{h}(u, p, q, m, \tilde{u}, \tilde{p}, \tilde{q}, \tilde{m}) &:= (G(m), \tilde{m}) \nonumber \\
&+ \int_{0}^{T} \int_{\Omega} \left( \frac{\partial u}{\partial t} \tilde{p} + \kappa \nabla u \cdot \nabla \tilde{p} + (\vec{v} \cdot \nabla u) \tilde{p} \right) dxdt \nonumber \\
&+ \int_{\Omega} \tilde{q}(x) \left( u(0, x) - m^{2}(x) \right) dx \nonumber \\
&+ \frac{1}{\sigma^{2}}\sum_{\ell=1}^{n}\sum_{k \in \mathcal{K}_{\ell}}\int_{0}^{T} \left( \mathcal{B}_\ell u - \mathbf{y}^{(k)} \right)^{T} \left( \mathcal{B}_\ell\tilde{u} \right) \delta(t - t_k) dt \nonumber \\
& + \int_{0}^{T}\int_{\Omega} \left( - \tilde{u} \frac{\partial p}{\partial t} + \kappa \nabla \tilde{u} \cdot \nabla p + (\vec{v} \cdot \nabla \tilde{u}) p \right) dxdt 
\end{align}
where $(u, p, q, m)$ are the variables from the gradient evaluation, and $(\tilde{u}, \tilde{p}, \tilde{q}, \tilde{m}) \in \mathcal{U} \times \mathcal{P} \times H^{1}(\Omega) \times \mathcal{M}$. Here, $\tilde{m}$ is a given direction, and we enforce $\tilde{q}(x) = \tilde{p}(0, x)$ in $\Omega$.

We require the first variation of $\mathcal{L}^h$ with respect to $p$ to vanish for all test functions $\hat{p} \in \mathcal{P}$:
\begin{align}\label{eq: incremental forward problem weak form}
0 &= \delta_{p} \mathcal{L}^{h} \hat{p} = -\int_{\Omega} 2 \hat{p}(0, x) m(x) \tilde{m}(x) dx + \int_{0}^{T} \int_{\Omega} \left( \frac{\partial \tilde{u}}{\partial t} \hat{p} + \kappa \nabla \tilde{u} \cdot \nabla \hat{p} + (\vec{v} \cdot \nabla \tilde{u}) \hat{p} \right) dxdt,
\end{align}
where we used integration by part and set the initial condition $\tilde{u}(0, x) = 0$.
This equation constitutes the weak form of the incremental forward problem, which we solve to obtain the incremental state $\tilde{u}$.

Next, we require the first variation of $\mathcal{L}^h$ with respect to $u$ to vanish for all test functions $\hat{u} \in \mathcal{U}$:
\begin{align}
0 = \delta_{u} \mathcal{L}^{h} \hat{u} = &\int_{0}^{T} \int_{\Omega} \left( -\frac{\partial \tilde{p}}{\partial t} \hat{u} + \kappa \nabla \hat{u} \cdot \nabla \tilde{p} + (\vec{v} \cdot \nabla \hat{u}) \tilde{p} \right) dxdt \nonumber \\
& + \frac{1}{\sigma^{2}}\sum_{\ell=1}^{n}\sum_{k \in \mathcal{K}_{\ell}}\int_{0}^{T} \left( \mathcal{B}_\ell \hat{u} \right)^{T} \left( \mathcal{B}_\ell\tilde{u} \right) \delta(t - t_k) dt,
\end{align}
where we used integration by part and set the terminal condition $\tilde{p}(T, x) = 0$ and $\tilde{q}(x) = \tilde{p}(0, x)$.
This equation constitutes the weak form of the incremental adjoint problem, which determines the incremental adjoint variables $\tilde{p}$ by solving the incremental adjoint problem backward in time.

Finally, with the computed states and adjoints, we derive the Hessian action. The weak form of the Hessian of $\mathcal{J}(m)$ evaluated at $m$ and acting in the direction $\tilde{m}$ is given by:
\begin{align}
(H(m) \tilde{m}, \hat{m}) = \delta_{m} \mathcal{L}^{h} \hat{m} = &\int_{\Omega} \mathcal{A}(\hat{m})(x) \mathcal{A}(\tilde{m})(x) dx \nonumber \\
& - \int_{\Omega} p(0, x) 2\hat{m}(x) \tilde{m}(x) dx - \int_{\Omega} \tilde{p}(0, x) 2m(x) \hat{m}(x) dx.
\end{align}
Eliminating the arbitrary test function $\hat{m}$, we obtain the strong form of the Hessian action:
\begin{equation}\label{eq: cost functional Hessian evaluation}
H(m) \tilde{m} =
\begin{cases}
\mathcal{A}^{2}(\tilde{m}) - 2p(0, \cdot)\tilde{m} - 2\tilde{p}(0, \cdot)m, & \text{in } \Omega, \\
\gamma \nabla \tilde{m} \cdot \vec{n} + \beta \tilde{m}, & \text{on } \partial \Omega.
\end{cases}
\end{equation}
For Gaussian--Newton approximation of the Hessian as used in \eqref{eq: GN approximation} for the Laplace approximation of the posterior covariance, we omit the second-order derivative term $- 2p(0, \cdot)\tilde{m}$ in \eqref{eq: cost functional Hessian evaluation}.

\subsection{Computation of the projected Jacobian and active subspace of the PDE solution} \label{subsection: Jacobian computation for neural operator}

As discussed in Sections~\ref{subsection: Dimension reduction} and~\ref{subsection: Data generation and derivative-informed training}, computing the expected cumulative Jacobian in \eqref{eq: active subspace} and the projected Jacobians in \eqref{eq: data projection} is critical for two purposes: performing active subspace dimension reduction and generating derivative-informed training datasets for the LANO architecture. We provide the mathematical derivations of these quantities for the model \eqref{eq: time-dependent advection-diffusion equation} below.

Given a random field parameter $m$, we derive the sensitivity equation by differentiating the governing PDE and initial condition \eqref{eq: time-dependent advection-diffusion equation} with respect to $m$ in an arbitrary direction $\hat{m}$. Let $\hat{u}$ denote the sensitivity of the state variable (the Fréchet derivative of $u$ in the direction $\hat{m}$). Since the parameter $m$ enters the model solely through the initial condition, the linearized initial condition is given by:
\begin{equation}
\label{eq: Jacobian computation initial condition derivative}
\hat{u}(0) :=  \nabla_m u(0) \hat{m} = 2m\hat{m}, \quad \text{in } \Omega.
\end{equation}
Moreover, because the differential operator in the PDE is independent of $m$, the linearized PDE governing $\hat{u}$ takes the same form as the state problem \eqref{eq: time-dependent advection-diffusion equation}. The system for Jacobian action $\hat{u}(t):=\nabla_m u(t) \hat{m}$ is:
\begin{align}
    \frac{\partial \hat{u}}{\partial t} - \kappa\Delta \hat{u} + \vec{v}\cdot \nabla \hat{u} = 0, & \text{ in } \Omega \times (0, T), \nonumber \\
    \label{eq: Jacobian computation sensitivity PDE}
    \hat{u}(0) = 2m\hat{m}, & \text{ in } \Omega, \\
    \nabla \hat{u} \cdot \vec{n} = 0, & \text{ on } \partial\Omega \times (0, T). \nonumber
\end{align}
Hence, the training data for the projected Jacobian in \eqref{eq: data projection} is computed in two steps: (1) we calculate the Jacobian action $\nabla_m u(t) \hat{m}$ in directions $\hat{m} = \psi_m^j$ for $j = 1, \dots, r_m$, where $\{\psi_m^j\}$ are the eigenvectors obtained from the eigenvalue problem \eqref{eq: active subspace}; and (2) we perform a left projection of these actions using $\boldsymbol{\Psi}_u^\top$.

To solve the eigenvalue problem \eqref{eq: active subspace} using randomized algorithms \cite{villa2021hippylib}, we must compute $(\nabla_m u(t))^* \hat{u} = (\nabla_m u(t))^* (\nabla_m u(t)) \hat{m}$ for random direction $\hat{m} \in L^2(\Omega)$. Here, $(\nabla_m u(t))^*$ denotes the adjoint of the Jacobian operator $\nabla_m u(t)$, defined by the relationship:
\begin{equation}
\label{eq:adjoint-relation}
\langle (\nabla_m u(t))^* v, \hat{m} \rangle = \langle v, \nabla_m u(t) \hat{m} \rangle = \langle v, \hat{u}(t) \rangle, \quad \forall v, \hat{m} \in L^2(\Omega),
\end{equation}
where $\langle \cdot, \cdot \rangle$ represents the $L^2(\Omega)$ inner product.
Let $p \in \mathcal{P}$ be an adjoint variable. From the sensitivity equation \eqref{eq: Jacobian computation sensitivity PDE}, we have:
\begin{equation}
    \int_0^t \int_\Omega \left(\frac{\partial \hat{u}}{\partial \tau} - \kappa\Delta \hat{u} + \vec{v}\cdot \nabla \hat{u}\right) p \, dx \, d\tau = 0.
\end{equation}
Applying integration by parts twice and using the properties $\nabla \cdot \vec{v} = 0$ in $\Omega$ and $\vec{v} \cdot \vec{n} = 0$ on $\partial \Omega$ by definition of the velocity field, we obtain:
\begin{equation}
\int_{\Omega} \hat{u}(t) p(t) \, dx = \int_{\Omega} \hat{u}(0) p(0) \, dx,
\end{equation}
provided that $p$ satisfies the following adjoint equation evolving backward in time:
\begin{align}
    - \frac{\partial p}{\partial \tau} - \kappa\Delta p - (\vec{v}\cdot \nabla) p = 0, & \text{ in } \Omega \times (0, t), \nonumber \\
    \label{eq: Jacobian computation adjoint PDE}
    p(t) = \hat{u}(t), & \text{ in } \Omega, \\
    \nabla p \cdot \vec{n} = 0, & \text{ on } \partial\Omega \times (0, t). \nonumber
\end{align}
By setting $v = \hat{u}(t)$ in \eqref{eq:adjoint-relation}, we derive:
\begin{equation}
\langle (\nabla_m u(t))^* \hat{u}(t), \hat{m} \rangle = \langle \hat{u}(t), \hat{u}(t) \rangle = \langle \hat{u}(0), p(0) \rangle = \langle 2m p(0), \hat{m} \rangle,
\end{equation}
where the last equality follows from the linearized initial condition \eqref{eq: Jacobian computation initial condition derivative}. This yields the final result: 
\begin{equation}
    (\nabla_m u(t))^* \hat{u}(t) = 2m p(0).
\end{equation}


\section{Numerical experiments} \label{section: Numerical experiments}

This section reports numerical results for the contaminant source tracking case study (Section~\ref{section: Case study: contaminant source tracking}). We (i) construct and validate the derivative-informed LANO surrogate, including dimension reduction, training setup, and accuracy for both states and Jacobians; (ii) benchmark surrogate-accelerated Bayesian inversion against a high-fidelity FEM baseline in terms of D-optimality quality and cost, and assess the efficient D-optimality evaluation that uses prior samples as MAP proxies; and (iii) apply the surrogate-accelerated PGRL pipeline (Algorithm~\ref{alg: pgrl}) to learn and interpret adaptive policies for sequential sensor placement.

\subsection{Surrogate modeling}

The spatial mesh, visualized in Figure~\ref{fig: physical domain + velocity fields a}, undergoes two levels of refinement. We employ a space of continuous, piecewise linear functions to discretize both the state variable and the random field parameter. This results in a high-dimensional discretization with dimension 7,863.

\begin{figure}[!htb]
    \centering
    \begin{subfigure}[t]{0.45\linewidth}
        \centering
        \includegraphics[width=\linewidth]{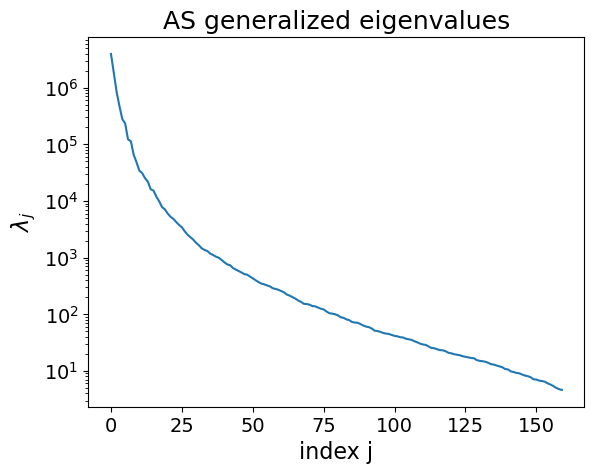}
    \end{subfigure}
    \hspace{0.25em}
    \begin{subfigure}[t]{0.45\linewidth}
        \centering
        \includegraphics[width=\linewidth]{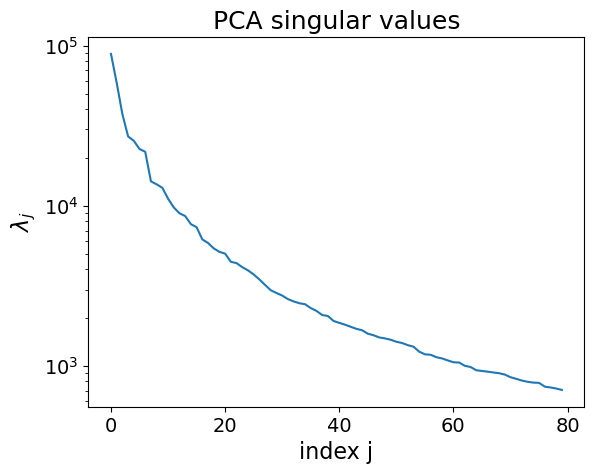}
    \end{subfigure}
    \caption{Decay of the generalized eigenvalues of active subspace projection for random field parameter dimension reduction (left) and the singular values of PCA for state variable dimension reduction (right).}
    \label{fig: decay of singular values and generalized eigenvalues}
\end{figure}

To reduce the dimension of the discretized random field parameter via active subspace projection, we sample $100$ parameter realizations from the prior to estimate the expectation in \eqref{eq: active subspace}. The Jacobian computation follows the derivations in Section~\ref{subsection: Jacobian computation for neural operator}. This approximate generalized eigenvalue problem is solved using the double-pass randomized algorithm \cite{villa2021hippylib}, with the decay of the $r_m = 160$ largest generalized eigenvalues shown in Figure~\ref{fig: decay of singular values and generalized eigenvalues}, leading to a parameter projection error of $2.48\%$.

For PCA projection of the state variables, we sample $M_P = 2000$ realizations of the random field parameter from the prior and compute the snapshot matrix \eqref{eq: snapshot matrix}. For the truncated SVD \eqref{eq: SVD}, we select the largest $r_u = 80$ singular values, which account for over 99.8\% of the total variance in the centered snapshot matrix $\hat{\mathbf{U}}$ and $1.24\%$ state projection error. The decay of these singular values is visualized in Figure~\ref{fig: decay of singular values and generalized eigenvalues}.

In the latent encoding layer, the dimension of the latent variables is set to $d_z = 200$, using a $\tanh$ activation function ($\sigma_z$). The latent attention mechanism employs projection matrices of dimension $d_a = 100$. For the latent dynamics, we use the ELU activation function ($\sigma_{\mathcal{L}}$). Finally, for the reconstruction of the reduced state variables, the hidden dimension is set to $r_u' = 128$.

For dataset generation, we reuse the $M_{P} = 2000$ parameter samples and state variable solutions obtained in dimension reduction. The projected Jacobian computation follows the derivations in Section~\ref{subsection: Jacobian computation for neural operator}. The training dataset consists of $M_T = 1600$ trajectories from the $M_P = 2000$ generated samples. The remaining samples are equally divided into validation and testing sets. The balancing hyperparameter $\lambda$ in the training loss \eqref{eq: training loss} is set to 1.0.

The LANO model is trained for 800 epochs with a batch size of 100. We utilize forward-mode automatic differentiation to compute the Jacobian with respect to the input reduced parameter. The model is evaluated on the validation dataset after each epoch. A learning rate scheduler reduces the learning rate by a factor of 0.5 if the validation loss plateaus for 10 consecutive epochs.

\begin{figure}[!htb]
  \centering
  
  \begin{subfigure}[t]{0.45\linewidth}
    \centering
    \includegraphics[width=\linewidth]{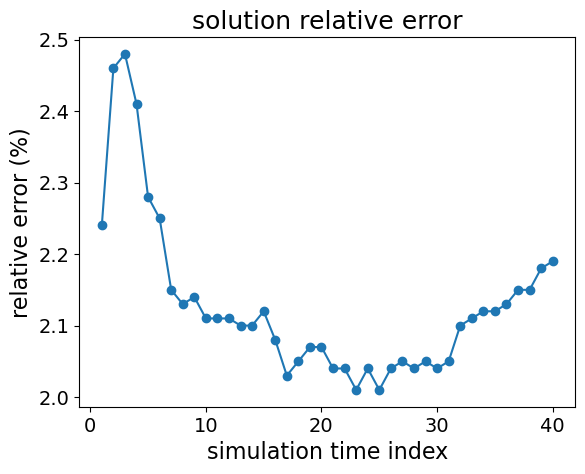}
  \end{subfigure}
  \hspace{0.25em}
  \begin{subfigure}[t]{0.45\linewidth}
    \centering
    \includegraphics[width=\linewidth]{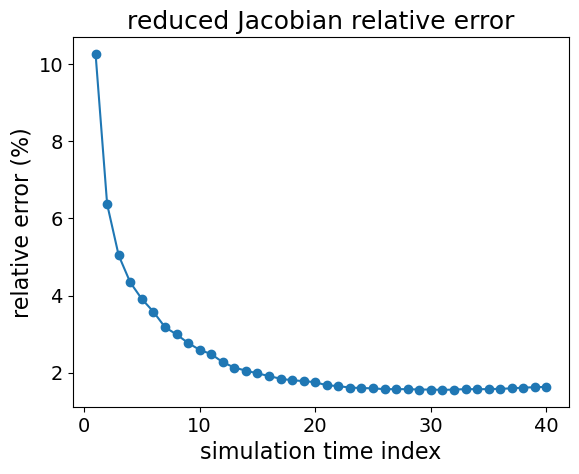}
  \end{subfigure}

    \caption{Relative errors of trained LANO on the testing dataset. The reduced Jacobian relative error (right) peaks at the first simulation time and gradually decreases. The solution relative error (left) shows less variation and no monotonic trend.}
  \label{fig: LANO testing}
\end{figure}

Post-training, the LANO model is evaluated on the testing dataset. Figure~\ref{fig: LANO testing} reports the average relative errors for both the solution and the reduced Jacobian, averaged over the simulation time at 2.13\% and 2.42\%, respectively. The high initial relative error in the reduced Jacobian is primarily due to the parameter $m$ entering the system exclusively through the initial condition $u(0, x) = m^2(x)$. At $t \approx 0$, the Jacobian is highly sensitive to the spatial details that are difficult for the surrogate to approximate accurately.

\begin{figure}[!htb]
  \centering
  \includegraphics[width=\linewidth]{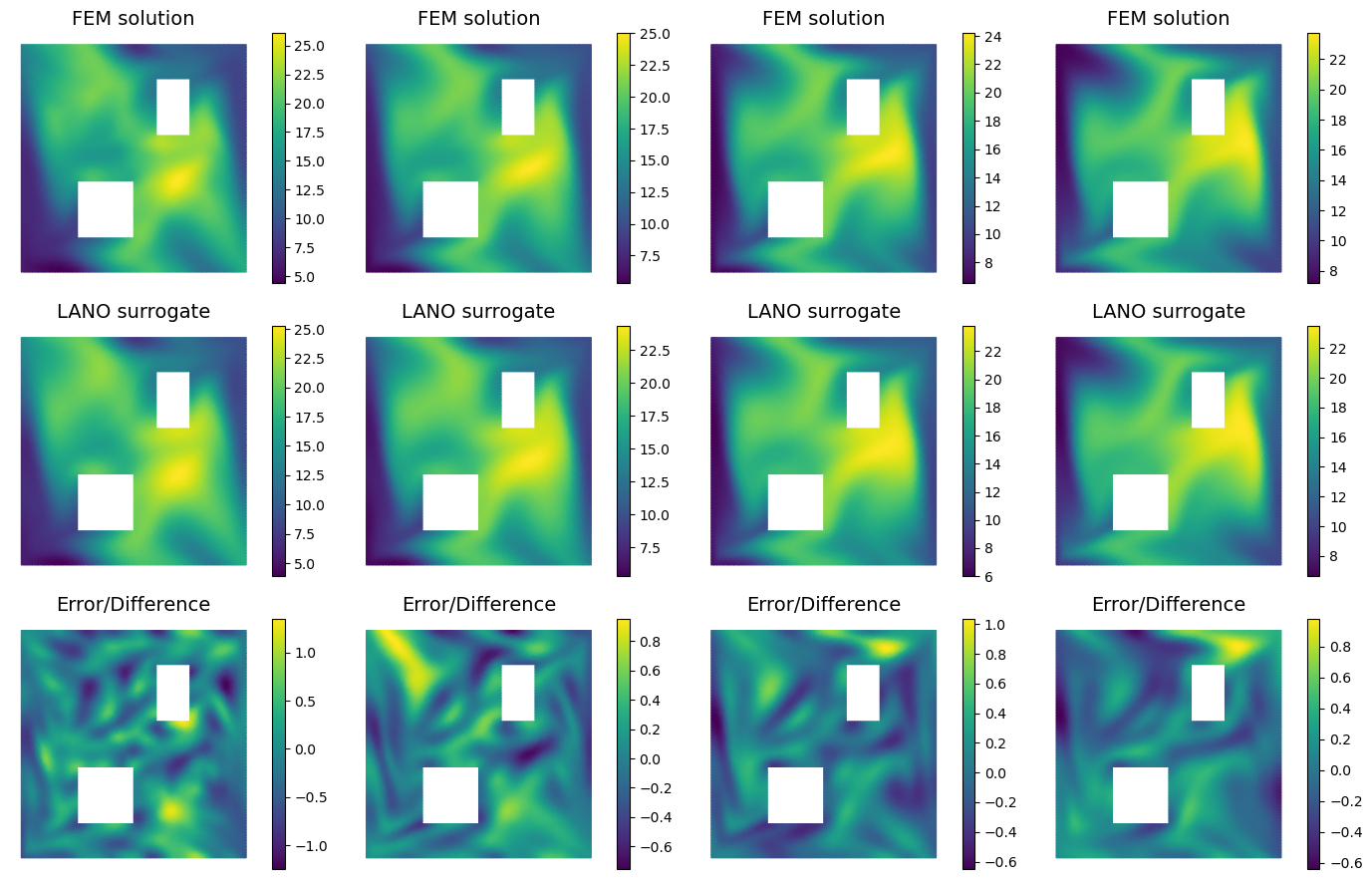}
  \caption{Contaminant concentration $u$ at simulation times 0.4, 0.8, 1.2, and 1.6 (from left to right). The LANO surrogate (middle) closely matches the FEM reference (top). Errors (bottom) indicate LANO effectively captures the dynamics.}
  \label{fig: FEM-NN-error}
\end{figure}

Figure \ref{fig: FEM-NN-error} provides a visual comparison between the high-fidelity FEM simulations and the LANO surrogate predictions for the contaminant concentration $u$ at time $t = 0.4, 0.8, 1.2, 1.6$. The LANO model (middle row) accurately reproduces the spatiotemporal evolution of the contaminant observed in the FEM reference (top row) across all four simulation time steps. The corresponding residual plots (bottom row) demonstrate that the approximation error remains small throughout the simulation, confirming that the surrogate effectively captures the underlying transport dynamics.

\subsection{Efficient computation of the D-optimality}

To compute the D-optimality \eqref{eq: terminal formulation terminal} in the terminal formulation of the reward function \eqref{eq: SBOED MDP formulation}, we use the approximation \eqref{eq: D-optimality computation} based on the Laplace approximation of the posterior covariance at the MAP point \eqref{eq: posterior covariance exact}, where the MAP point is computed by solving the optimization problem \eqref{eq: MAP optimization} and the eigenvalues are computed by solving the generalized eigenvalue problem \eqref{eq: generalized eigenvalue}. We employ a Newton-CG FEM solver initialized at the prior mean to compute the MAP point using hiPPYlib \cite{villa2021hippylib}, with the gradient and Hessian defined in Section \ref{sec:gradHess}. The solver is configured with relative and absolute tolerances of $10^{-12}$ and $10^{-6}$, respectively, and a maximum of 50 Newton iterations. Subsequently, the generalized eigenvalue problem \eqref{eq: generalized eigenvalue} is solved via a randomized double-pass algorithm with target rank $2K'$ and an oversampling parameter of 5, where $K'$ is the number of simulation time steps used in the forward model. We also solve the projected optimization problem \eqref{eq: accelerated MAP} using the trained LANO surrogate. The optimization variable $\boldsymbol{\beta}_m$ is initialized as a zero vector. We employ the L-BFGS optimizer, configured with an initial step size of $0.5$, a maximum of $200$ iterations, a history size of $200$, and a strong Wolfe line search. For the eigenvalue computation, we solve the projected eigenvalue problem \eqref{eq: reduced geneig} with LANO surrogate of the Jacobian.

\begin{figure}[htbp]
    \centering
        \includegraphics[width=0.48\linewidth]{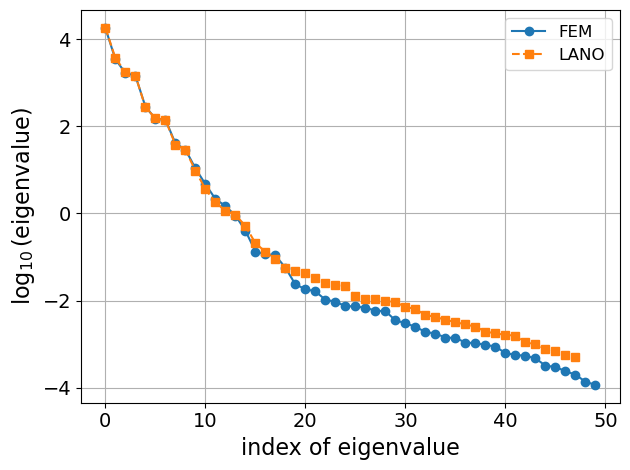}
        \includegraphics[width=0.48\linewidth]{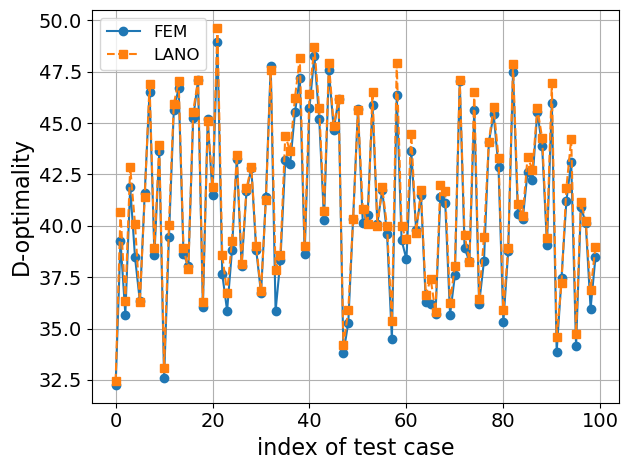}
    \caption{Left: Comparison of the eigenvalues of the Gauss--Newton Hessian in \eqref{eq: generalized eigenvalue} computed by FEM  and in \eqref{eq: reduced geneig} computed by LANO surrogate at one random test sample. Right: Comparison of the D-optimality \eqref{eq: D-optimality computation} at 100 random test samples.}
    \label{fig: KL D-optimality comparison}
\end{figure}

To access the accuracy of the surrogate-accelerated D-optimality computation compared to the FEM reference, we use $500$ total test cases, generated by sampling $5$ random field parameters from the prior and pairing each with $100$ distinct sensor placements obtained by randomly perturbing the initial configuration (Figure~\ref{fig: initial placement of sensors}). The relative root mean square error (RMSE) of the D-optimality values across all test cases is 1.29\% with a standard deviation of the relative error at 1.25\%, which is satisfactory for design purposes. The comparison of the eigenvalues of the Gauss--Newton Hessian in \eqref{eq: generalized eigenvalue} computed by FEM and in \eqref{eq: reduced geneig} computed by LANO surrogate at one random test sample is shown in the left panel of Figure~\ref{fig: KL D-optimality comparison}. The right panel of Figure~\ref{fig: KL D-optimality comparison} presents the comparison of the D-optimality values at 100 random test samples, demonstrating close agreement between the LANO approximation and the high-fidelity reference.

\begin{figure}[!htb]
  \centering
  \includegraphics[width=0.9\linewidth]{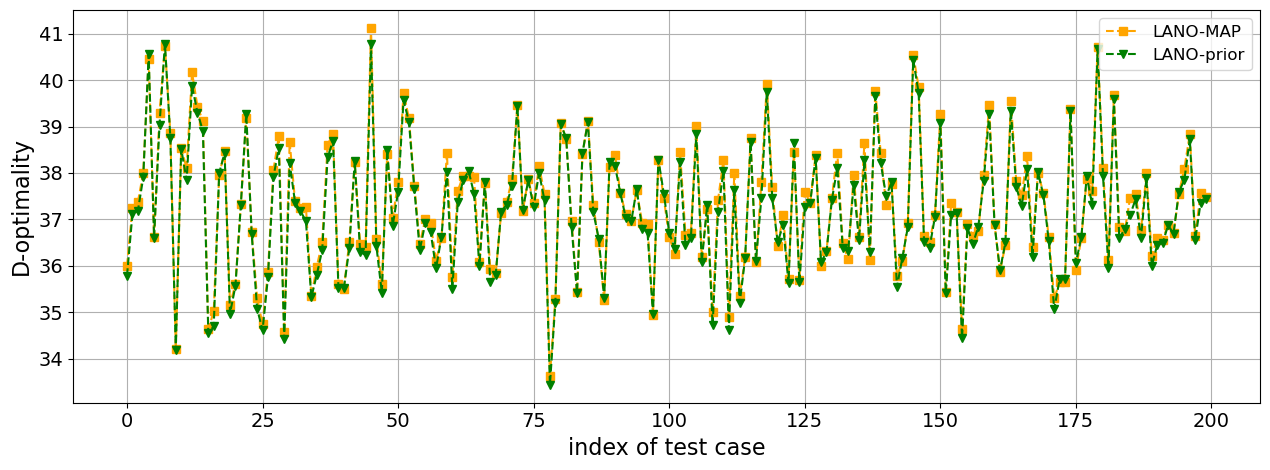}
    \caption{Comparison of D-optimality computed by LANO using the standard approach (solving for MAP first) and the efficient strategy (using prior samples), over 200 test cases for the first parameter, with $0.44\%$ relative RMSE and $0.99$ correction.}
  \label{fig: D-optimality comparison efficient}
\end{figure}

We note that the MAP point acceleration by LANO is not significant, 8.93 seconds by FEM versus 5.65 seconds by LANO, since the optimization still requires multiple iterations. The main acceleration arises from the efficient eigenvalue computation using LANO, which takes only 0.14 seconds compared to 13.46 seconds using FEM, about $100\times$ speedup, as it avoids solving large-scale PDEs.
To further reduce the cost of computing reduced MAP points, we adopt an efficient D-optimality evaluation strategy: we compute the eigenvalues at each prior sample and treat that sample as a proxy for the MAP point (for the synthetic data generated from it), thereby avoiding repeated MAP optimizations inside the design loop. To validate this strategy, we drew $5$ random field parameters from the prior and generated $200$ distinct sensor placements by perturbing the initial configuration $N=4$ times within the design bounds. For each test case, we computed D-optimality using (i) the standard approach (solve for the MAP, then compute eigenvalues) and (ii) the efficient approach (use the prior sample directly). The standard approach serves as the reference. Figure~\ref{fig: D-optimality comparison efficient} visualizes the comparison for the first parameter. The relative RMSE and correlation across all test cases are $0.44\%$ and $0.99$, respectively, indicating that the efficient approximation of the D-optimality at the prior sample closely tracks the D-optimality computed at the corresponding MAP point.

\begin{figure}[!htb]
  \centering  
    \includegraphics[width=\linewidth]{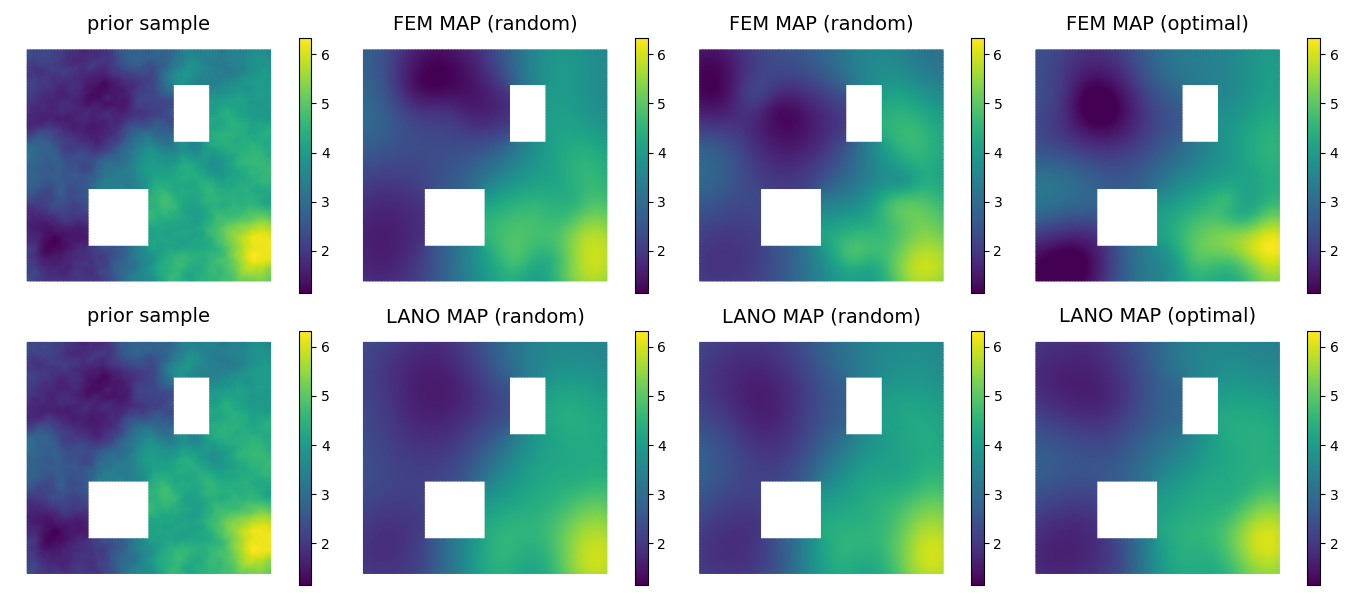}
  \caption{Visual comparison of a prior sample (left) against MAP estimates for two random designs (middle) and the optimal design (right). Top row: High-fidelity FEM reference solutions. Bottom row: LANO surrogate approximations. The MAP estimates closely resemble the prior, with surrogate relative $L_2$-errors of $8.64\%$, $10.11\%$, and $10.22\%$, respectively.}
  \label{fig: MAP points random optimal}
\end{figure}

To further justify this approximation, we visually inspect the parameter fields in Figure~\ref{fig: MAP points random optimal}. The figure displays a representative prior sample alongside the MAP estimates obtained for two random sensor configurations and the optimal design. Qualitatively, the MAP estimates exhibit spatial structures that are similar to the underlying prior sample. This proximity in the parameter space implies that the local curvature of the log-posterior, and consequently the spectrum of the Hessian, does not vary significantly between the prior sample and the posterior mode. This observation explains the high accuracy of the proxy strategy reported above. Furthermore, the LANO surrogate (bottom row) faithfully reproduces the high-fidelity FEM reference solutions (top row) with relative $L_2$ errors between $8.64\%$ and $10.22\%$, confirming its reliability for these accelerated design evaluations.

\subsection{SBOED via PGRL for optimal sensor placement}

We summarize the network architectures, data generation procedure, and optimization hyperparameters used to train the policy and Q-networks in Algorithm~\ref{alg: pgrl}.

\emph{Network Architecture:}
We employ MLPs for both the policy and Q-networks. The policy network comprises two hidden layers, each with a width set to ten times the input dimension, utilizing ReLU activation functions. To ensure the design vectors respect the prescribed bounds, the final output is passed through a sigmoid function. The Q-network mirrors this architecture, utilizing the same depth and activation functions, but differs in its input/output dimensions and omits the final constraint on the output.

\emph{Data Generation:}
To accelerate training, the data collection component of Algorithm~\ref{alg: pgrl} leverages offline precomputation. Since the prior distribution and the trained LANO surrogate remain fixed, we pre-sample random field parameters and compute the corresponding state trajectories and Jacobian approximations offline. These quantities are organized into $L=200$ batches, each containing $M=1000$ trajectories. During online training, one batch is loaded per policy update step to efficiently evaluate the design vectors and Hessians required for D-optimality.

\emph{Training Strategy:}
Optimization is performed using the Adam optimizer \cite{kingma2017adam}. The policy network is trained with an initial learning rate of $1 \times 10^{-4}$ over $L=200$ global updates. The Q-network employs a learning rate of $3 \times 10^{-4}$ and undergoes $L'=200$ internal updates for every single policy step. Both networks utilize an exponential scheduler, decaying the learning rate by a factor of $0.98$ at each step.

With the configuration described above, the total execution time is $7,261$ seconds, averaging approximately $36$ seconds per policy update. The primary computational bottleneck is the construction of pointwise observation operators, which must be regenerated whenever sensor coordinates change. For each of the $M=1000$ trajectories in a batch, $N=4$ operators are generated, resulting in $4,000$ sequential constructions per step. This process consumes approximately $30$ seconds of the $36$-second update cycle.

\begin{figure}[!htb]
  \centering
    \begin{subfigure}[t]{0.31\linewidth}
        \centering
        \includegraphics[width=\linewidth]{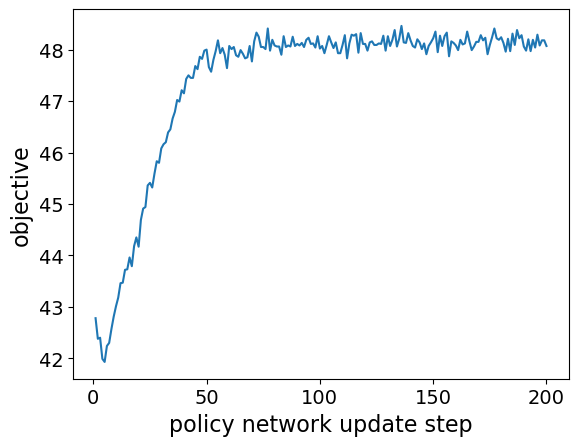}
    \end{subfigure}
    \begin{subfigure}[t]{0.31\linewidth}
        \centering
        \includegraphics[width=\linewidth]{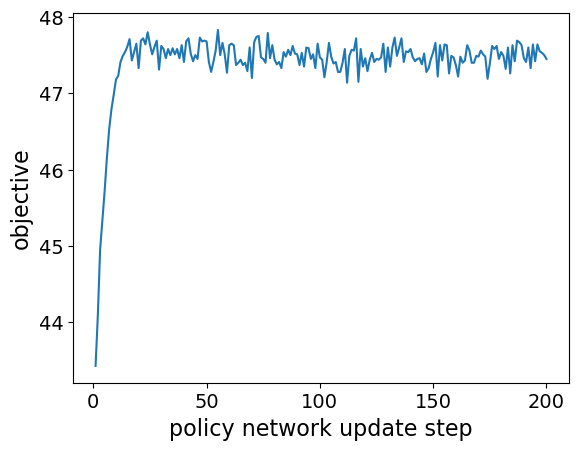}
    \end{subfigure}
    \begin{subfigure}[t]{0.31\linewidth}
        \centering
        \includegraphics[width=\linewidth]{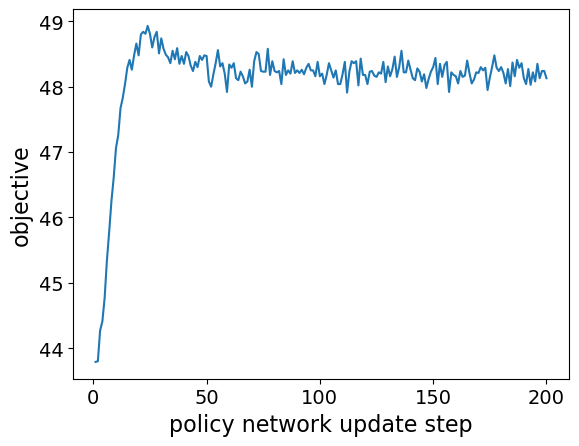}
    \end{subfigure}
    \caption{Evolution of the objective value with respect to the policy network update step for three different random seeds.}
    \label{fig: objective evolution} 
\end{figure}

Figure~\ref{fig: objective evolution} illustrates the evolution of the objective value over the course of training with three different random seeds for initialization. The rapid initial increase indicates swift policy improvement, after which the objective stabilizes. The minor fluctuations observed in the later stages are characteristic of the stochastic nature of the gradient updates and do not indicate divergence.

We evaluate the efficacy of the optimized policy (Algorithm~\ref{alg: pgrl}) by applying it to $5$ distinct random field parameters randomly sampled from the prior. For each parameter, the trained policy network adaptively adjusts the sensor placements, and the resulting D-optimality is computed. To benchmark performance, we compare these policy-driven designs against a baseline of random sensor placements. Specifically, for each parameter, we generate $200$ random sensor configurations using the same perturbation strategy employed in the efficient D-optimality validation. Table~\ref{tab: Comparison to random sensor displacement} summarizes the results, with the third column indicating the percentage of random configurations outperformed by the learned policy.

\begin{table}[!htb]
    \centering
    \caption{Comparison of policy-based sensor placement against 200 random displacements.}
    \label{tab: Comparison to random sensor displacement}
    \begin{tabular}{ccc}
        \toprule
        \textbf{Parameter Index} & \textbf{D-optimality} ($\uparrow$) & \textbf{Win Rate} (\%) ($\uparrow$) \\
        \midrule
        1 & 40.39 & 97.5 \\
        2 & 52.97 & 98.0 \\
        3 & 46.42 & 97.5 \\
        4 & 46.56 & 98.0 \\
        5 & 46.75 & 98.0 \\
        \bottomrule
    \end{tabular}
\end{table}

\begin{figure}[!htb]
  \centering
  
  \begin{subfigure}[t]{0.31\linewidth}
    \centering
    \includegraphics[width=\linewidth]{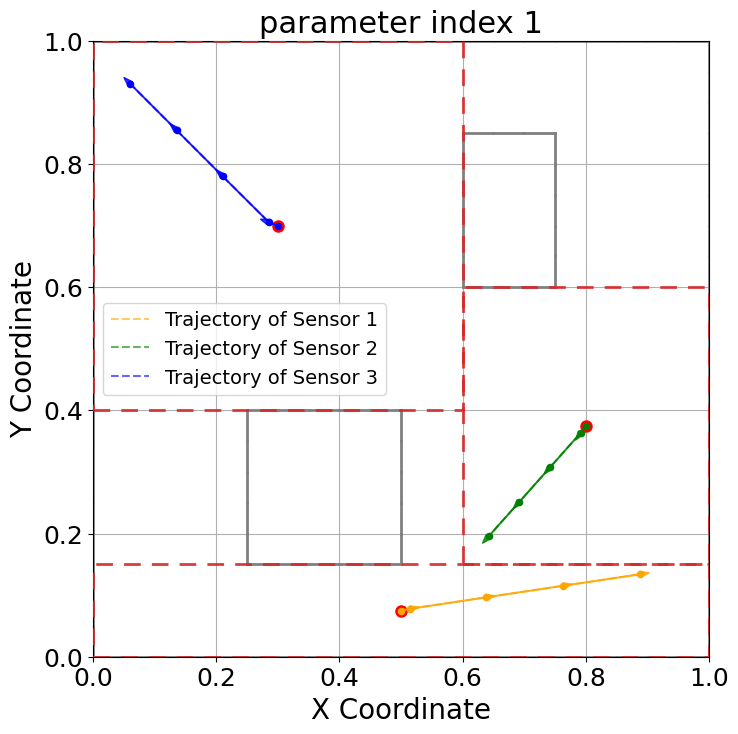}
  \end{subfigure}
  \hspace{0.25em}
  \begin{subfigure}[t]{0.31\linewidth}
    \centering
    \includegraphics[width=\linewidth]{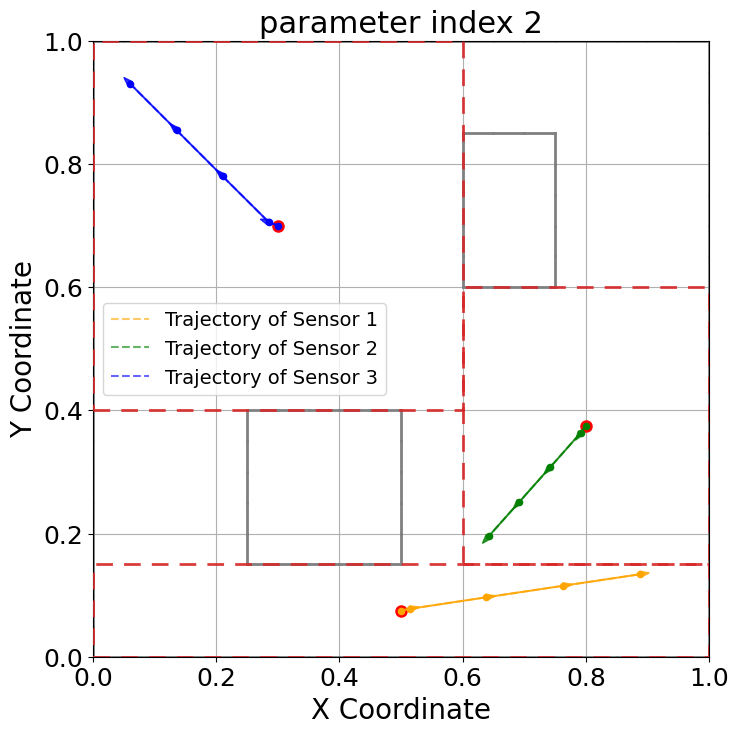}
  \end{subfigure}
  \hspace{0.25em}
  \begin{subfigure}[t]{0.31\linewidth}
    \centering
    \includegraphics[width=\linewidth]{figures/case-study/velocity_field_g1.png}
  \end{subfigure}
    
\caption{Visualization of optimized design vectors for two different parameter realizations (left and middle). The sensors consistently move upstream against the velocity field (right) in their admissible domains.}  
\label{fig: design vector g1}
\end{figure}

As shown in Table~\ref{tab: Comparison to random sensor displacement}, the policy-based designs consistently outperform the vast majority (over 97\%) of random placements, demonstrating the effectiveness of the learned strategy. Visualizing the resulting design vectors (Figure~\ref{fig: design vector g1}) reveals a consistent pattern across different parameter realizations. This consistency is expected, as the current framework trains a static policy over the prior distribution rather than adapting to real-time observational data. Consequently, the network learns a robust, generalized sensor placement strategy applicable to the entire class of random field parameters defined by the prior.

\begin{figure}[!htb]
  \centering
  
  \begin{subfigure}[t]{0.31\linewidth}
    \centering
    \includegraphics[width=\linewidth]{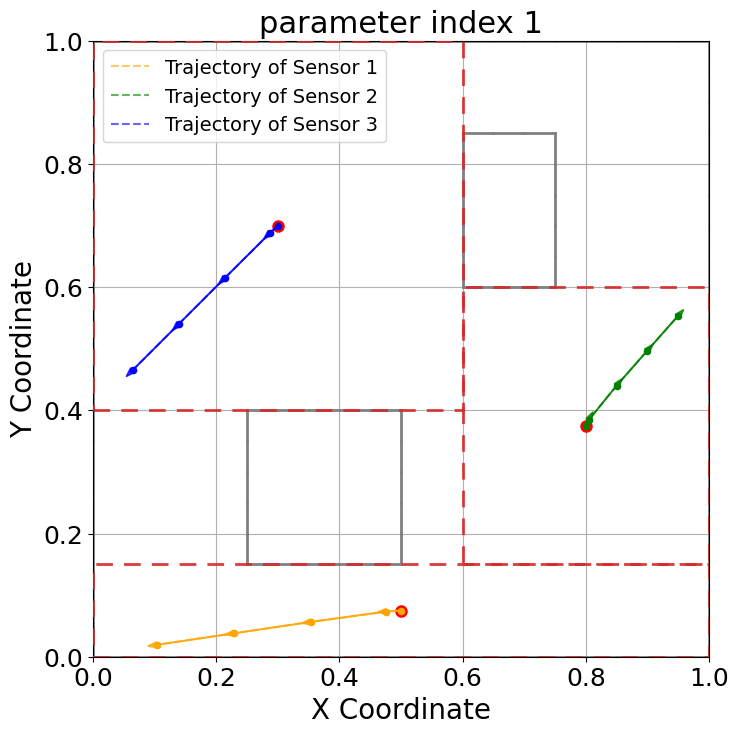}
  \end{subfigure}
  \hspace{0.25em}
  \begin{subfigure}[t]{0.31\linewidth}
    \centering
    \includegraphics[width=\linewidth]{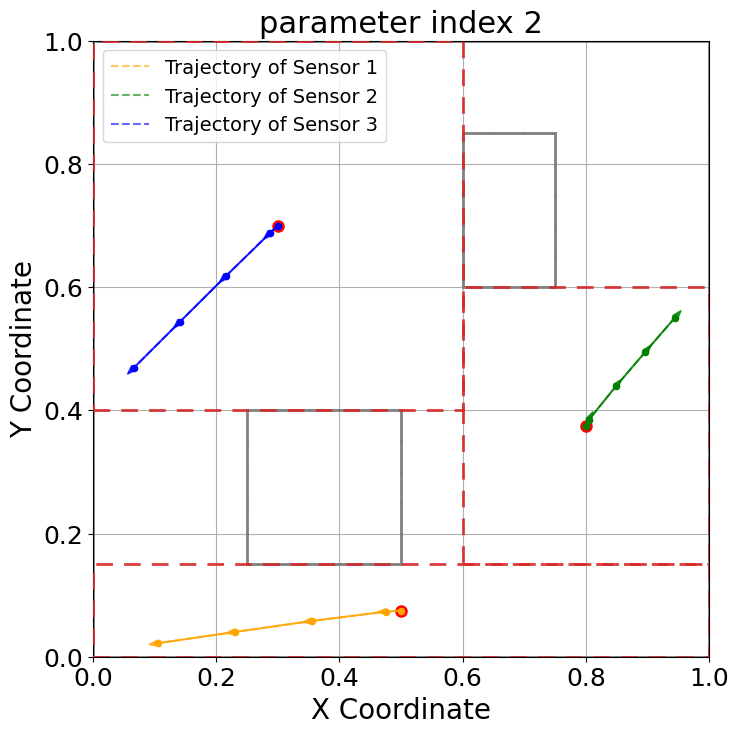}
  \end{subfigure}
  \hspace{0.25em}
  \begin{subfigure}[t]{0.31\linewidth}
    \centering
    \includegraphics[width=\linewidth]{figures/case-study/velocity_field_g2.png}
  \end{subfigure}
    
    \caption{Visualization of optimized design vectors for two different parameter realizations (left and middle) under a reversed velocity field (right). The policy correctly adapts to move sensors in the new upstream direction.}  
\label{fig: design vector g2}
\end{figure}

A physical interpretation of the learned strategy emerges when comparing the design vectors with the underlying velocity field (Figure~\ref{fig: design vector g1}). The sensors consistently align against the flow direction, effectively moving ``upstream.'' Intuitively, this maximizes information gain by positioning sensors closer to the contaminant source, thereby reducing diffusive blurring. To validate this hypothesis, we conducted a control experiment using a reversed velocity field (Figure~\ref{fig: physical domain + velocity fields c}). When Algorithm~\ref{alg: pgrl} is applied to this modified system, the policy adapts accordingly, shifting sensors in the new upstream direction (Figure~\ref{fig: design vector g2}). This adaptive behavior confirms that the policy has successfully encoded the underlying physics of the transport. 

\begin{figure}[!htb]
  \centering  
    \includegraphics[width=\linewidth]{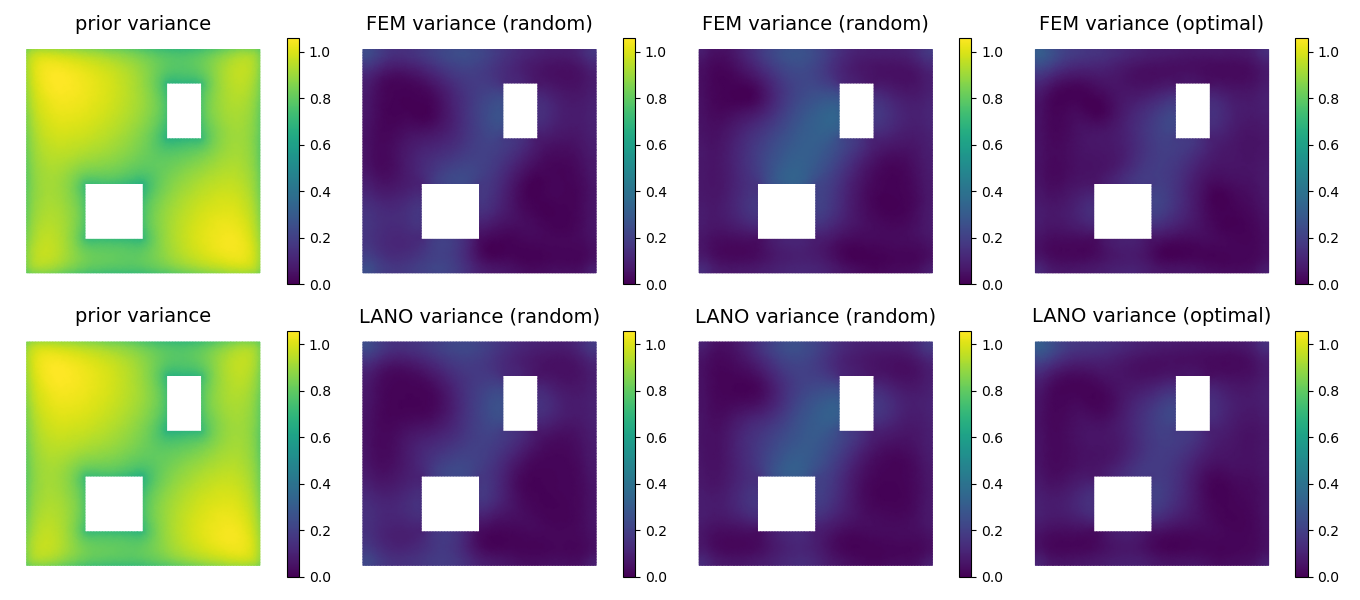}
  \caption{Comparison of posterior uncertainty reduction. From left to right: Pointwise variance of the prior, followed by the posterior variances for two random designs and the optimal design. The LANO surrogate (bottom row) closely reproduces the FEM reference variances (top row). The D-optimality estimates for the random and optimal designs are $49.15$, $46.28$, and $51.95$ (FEM), compared to $48.65$, $46.95$, and $52.19$ (LANO).}
  \label{fig: pointwise covariance}
\end{figure}

Finally, to visualize the impact of the learned policy on uncertainty reduction, we plot the pointwise variance fields in Figure~\ref{fig: pointwise covariance}. The leftmost column displays the prior variance, while the subsequent columns show the posterior variances for two random designs and the optimal design. While all sensor configurations reduce the initial uncertainty, the optimal design (rightmost column) yields a visibly greater reduction in variance across the domain compared to the random baselines. Moreover, the surrogate-predicted variance fields (bottom row) are visually indistinguishable from the high-fidelity FEM reference solutions (top row), providing further validation that the LANO surrogate accurately captures the posterior covariance structure.

\section{Conclusion} \label{section: conclusion}


We presented a fast and scalable framework for sequential Bayesian optimal experimental design (SBOED) in PDE-governed inverse problems with high-dimensional (potentially infinite-dimensional) random field parameters. The key idea is to amortize the design computation by casting SBOED as a finite-horizon Markov decision process and learning a policy offline via policy-gradient reinforcement learning (PGRL); once trained, the policy proposes new experiments online from the experiment history without repeatedly solving an SBOED optimization problem. To make both policy training and reward evaluation practical at scale, we combined dual dimension reduction (active subspace projection for parameters and principal component analysis for states) with an adjusted derivative-informed latent attention neural operator (LANO) surrogate. The surrogate predicts both the parameter-to-solution map and its Jacobian, enabling efficient sensitivity-aware computations. Using a Laplace approximation, we adopted a D-optimality reward and introduced an eigenvalue-based evaluation strategy that exploits low-rank structure and uses prior samples as proxies for maximum a posteriori (MAP) points, thereby avoiding expensive inner-loop MAP optimizations while retaining accurate information-gain estimates. 

Numerical experiments on contaminant source tracking validated the approach. The surrogate-accelerated workflow achieved approximately $100\times$ speedup over a high-fidelity finite element baseline while maintaining comparable posterior quality, and the learned sensor-placement policy consistently outperformed random designs (win rate $>97\%$). The resulting policies were also physically interpretable, discovering an ``upstream'' tracking strategy that counters diffusive blurring to maintain informativeness. 
Future work includes extending the approach to more complex nonlinear and multiphysics models, non-Gaussian priors and likelihoods, and robust objectives that account for model discrepancy; additional directions include incorporating hard design constraints and enabling online adaptation of the policy and surrogate under distribution shift.

\section*{Acknowledgement}

This work has been supported in part by NSF grants DMS-2245111 and CNS-2325631.

During the preparation of this work, the authors used ChatGPT and Google Gemini in order to polish the writing. After using this tool/service, the authors reviewed and edited the content as needed and take full responsibility for the content of the published article.

\bibliographystyle{elsarticle-num}
\bibliography{references}

\end{document}